\documentclass[twoside, 11pt]{article}

\usepackage{amssymb, amsmath, latexsym, mathrsfs, verbatim, calc}
\usepackage{color}

\numberwithin{equation}{section}
\newtheorem{thm}{Theorem}[section]
\newtheorem{lem}[thm]{Lemma}

\newtheorem{prop}[thm]{Proposition}
\newtheorem{rem}[thm]{Remark}

\newtheorem{defn}[thm]{Definition}
\newtheorem{assum}[thm]{Assumption}

\definecolor{blue-violet}{rgb}{0.54, 0.17, 0.89}
\definecolor{purple1}{rgb}{0.63, 0.36, 0.94}
\definecolor{purple2}{rgb}{0.87, 0.0, 1.0}
\definecolor{purple}{rgb}{0.5, 0.0, 0.5}
\definecolor{pansypurple}{rgb}{0.47, 0.09, 0.29}
\definecolor{orange}{rgb}{1.0, 0.27, 0.0}

\def\ba{\begin{array}}
\def\ea{\end{array}}
\def\beq{\begin{equation}}
\def\bes{\begin{equation*}}
\def\ees{\end{equation*}}
\def\bea{\begin{eqnarray}}
\def\eea{\end{eqnarray}}
\def\beaa{\begin{eqnarray*}}
\def\eeaa{\end{eqnarray*}}

\def\dis{\displaystyle}

\def\no{\noindent}

\def\lastline{\par \vspace{-7.3ex} \no}

\def\nts{\negthinspace}
\def\ss{\smallskip}
\def\ms{\medskip}
\def\bs{\bigskip}
\def\q{\quad}
\def\qq{\qquad}

\def\ol{\overline}

\def\ua{\mathop{\uparrow}}
\def\da{\mathop{\downarrow}}

\def\={=\nts \nts=\nts \nts=\nts \nts=}

\def\argmax{\mathop{\mbox{\rm argmax}}}

\def\({\textnormal{(}}
\def\){\textnormal{)}}

\def\cd{\cdot}
\def\cds{\cdots}

\def\pa{\partial}

\def\a{\alpha}
\def\b{\beta}
\def\g{\gamma}
\def\d{\delta}
\def\e{\varepsilon}

\def\k{\kappa}
\def\l{\lambda}
\def\m{\mu}

\def\si{\sigma}
\def\vsi{\varsigma}
\def\t{\tau}

\def\th{\theta}
\def\vth{\vartheta}
\def\o{\omega}

\def\vf{\varphi}
\def\p{\psi}

\def\vsi{\varsigma}
\def\D{\Delta}
\def\G{\Gamma}
\def\L{\Lambda}
\def\O{\Omega}

\def\P{\Psi}

\def\Th{\Theta}

\def\cF{{\cal F}}
\def\cG{{\cal G}}

\def\cO{{\cal O}}

\def\hC{\mathbb{C}}
\def\hD{\mathbb{D}}
\def\hE{\mathbb{E}}
\def\hF{\mathbb{F}}

\def\hN{\mathbb{N}}

\def\hP{\mathbb{P}}

\def\hR{\mathbb{R}}

\def\hX{\mathbb{X}}

\def\sB{\mathscr{B}}
\def\sC{\mathscr{C}}
\def\sD{\mathscr{D}}

\def\sF{\mathscr{F}}
\def\sG{\mathscr{G}}
\def\sH{\mathscr{H}}

\def\sL{\mathscr{L}}
\def\sM{\mathscr{M}}

\def\sU{\mathscr{U}}

\def\as{{\hbox{-a.s.}}}

\def\limsup{\mathop{\ol{\rm lim}}}

\def\neg{\negthinspace}

\def\no{\noindent}

\def\ss{\smallskip}
\def\ms{\medskip}
\def\bs{\bigskip}
\def\q{\quad}
\def\qq{\qquad}

\def\pa{\partial}
\def\cd{\cdot}
\def\cds{\cdots}

\def\qed{\hfill \rule[0cm]{.25cm}{.25cm}\medskip}   
\def\dfnn{\stackrel{\triangle}{=}}

\def\1{{\bf 1}}

\newenvironment{itm}{\vspace{-1ex}\begin{itemize}}{\end{itemize}}
\def\bi{\begin{itm}}
\def\ei{\end{itm}}

\def\equ_ind{\arabic{section}.\arabic{equation}}
\def\sec_ind{\arabic{section}}

\begin{document}
\title{\bf On Optimal Dividend and Investment Strategy  under Renewal Risk Models}
\author{Lihua Bai\thanks{\noindent School of Mathematics,
Nankai University, Tianjin, 300071, P. R. China. E-mail:
lbai@nankai.edu.cn. This author is supported in part by Chinese NSF
grants \#11471171 and \#11911530091. }, ~Jin Ma\thanks{\noindent
 Department of Mathematics,
University of Southern California, Los Angeles, CA 90089.
Email: jinma@usc.edu. This author is supported in part by NSF grant \#DMS-1908665. }\thanks{Part of this work was completed when both authors were visiting School of Statistics and Mathematics, Shanghai Lixin University of Accounting the Finance, Shanghai, China, whose hospitality was greatly appreciated.}}

\date{\today}
\maketitle

\vspace{-2mm}

\begin{abstract}
In this paper we continue investigating the optimal dividend and investment problems under the Sparre Andersen model. More precisely, we assume that the claim frequency is a renewal process instead of a standard compound Poisson process, whence semi-Markovian. Building on our previous work \cite{BaiMa17}, where we established the dynamic programming principle via a {\it backward Markovization} procedure and proved that the value
function is the unique {\it constrained} viscosity solution of the HJB equation, in this paper we focus on the construction of the optimal strategy. The
main difficulties in this effort is two fold: the regularity of the viscosity solution to a non-local, nonlinear, and degenerate parabolic PDE on an unbounded domain, which seems to be new in its own right; and the well-posedness of the closed-loop stochastic system. By introducing an
auxiliary PDE, we construct an $\e$-optimal strategy, and prove the well-posedness of the corresponding closed-loop system, via a ``bootstrap"
technique with the help of a Krylov estimate.
   \end{abstract}
\vfill \bs

\no

{\bf Keywords.} \rm Optimal dividend control, Sparrer-Anderson model, Hamilton-Jacobi-Bellman equation, viscosity solution, Krylov estimate.

\bs

\no{\it 2000 AMS Mathematics subject classification:} 60H07,15,30;
35R60, 34F05.

\eject

\section{Introduction}
In this paper we continue our investigation on the optimal dividend and investment problems under a Sparre Andersen insurance model. More
precisely, we assume that the claim number process is a {\it renewal} process instead of a standard Poisson process, therefore it is also referred to as a
 {\it renewal risk model}. Finding the optimal strategy for such a problem has been considered as an intriguing but challenging open problem for quite sometime (cf. e.g., \cite{AT} and reference cited therein)
mainly due to the semi-Markov nature of the renewal process, as well as the non-optimality of the well-known  barrier strategy
(see \cite{albrh}).  More specifically, for a general insurance model involving investments, even under the simplest Cram\'er-Lundberg form,
direct calculation of optimal strategy becomes almost impossible, and the solution procedure often depends on some
more general stochastic control technique. In particular, the approach of dynamic programming and consequently the study of the associated Hamilton-Jabobi-Bellman (HJB) equation and its viscosity solution, become a natural way to attack the problem (cf. e.g., \cite{azc, azcu}). However, as was pointed
out in \cite{AT}, the non-Markovian nature of Sparre Andersen model drastically complicated this approach, as it took away the basis of
dynamic programming.
On the other hand, since the commonly believed barrier type of dividend strategy was shown to be non-optimal in \cite{albrh}, the structure of the optimal dividend-investment strategy under a renewal risk model has naturally become an intriguing issue to explore.

Our recent paper \cite{BaiMa17} was the first step towards the solution to this problem. Specifically, in \cite{BaiMa17} we considered the following simplest (one-stock) ``toy"  model for the surplus process with dividend of the Sparre Andersen type.
Let $T>0$ be a given time horizon and $s\in[0,T]$, we consider the following SDE  on a filtered probability space $(\O, \cF, \hP; \hF)$):
 for $t\in [s,T]$,
%
\bea
\label{surplus0}
dX_{t}^{\pi}=pdt+rX_{t}^{\pi}dt+(\mu-r)\gamma_{t}X^\pi_{t}dt+ \sigma \gamma_{t}X_{t}^{\pi}dB_{t}-dQ_t-dL_{t}, \q X^\pi_s=x,
\eea
where $x$ is the initial surplus, $p$ the premium rate, $r$  the interest rate, and $(\mu, \si)$  the appreciation rate and the volatility of the stock,
respectively, all  assumed to be positive constants; $Q_t=\sum_{i=1}^{N_{t}}U_{i}$ is the (renewal) claim process, and $\pi=(\gamma_{t},L_{t})$, $t\ge0$,
is the investment-dividend pair, in which $\g=\{\g_t\}_{t\ge0}$ represents the proportion of the surplus invested in the stock at each time $t$ (hence $\gamma_{t}\in [0,1]$), and $L=\{L_t\}_{t\ge0}$ is the cumulative dividends process (hence increasing).
Denoting $\sU_{ad}$ to be all such investment-dividend strategies and the  solution to (\ref{surplus0}) by $X_t=X^\pi_t=X_{t}^{\pi,x}$,
define $\tau^{\pi}_s=\t^{\pi,x}_s:={\rm inf}\{t\geq s; X_{t}^{\pi,x}< 0\}$, $s\in[0,T]$, to be the ruin time of the insurance company. The goal is to maximize the following expected cumulated dividends:
\bea
\label{cost0}
J(s,x;\pi):=\hE\Big\{\int_{s}^{\tau^{\pi,x}_s\wedge T}e^{-c (t-s)}dL_{t}\Big\}:=\hE\Big\{\int_{s}^{\tau^{\pi}_s\wedge T}e^{-c (t-s)}dL_{t}\Big|X^\pi_{s}=x\Big\}, \q s\in[0,T],
\eea
where $c>0$ is the {\it discounting factor}.

We should note that although the problem (\ref{surplus0})-(\ref{cost0}) is the simplest model that contain both investment and dividend, the
solution of it is surprisingly challenging. The first obstacle is the fact that the claim frequency (or counting) process $N$ is renewal, hence non-Markovian (see \S2 for details). Thus the dynamic programming approach does not apply directly. To overcome this difficulty we follow a standard
``backward Markovization" procedure by adding an extra state variable $W=\{W_t\}_{t\ge0}$ that measures the time elapsed since the last claim (see \S2 for details) so that the model becomes Markovian, and the dynamic programming approach becomes valid. Along this line, in  \cite{BaiMa17} we verified the  {\it dynamic programming principle} (DPP), and proved that the value function of problem (\ref{surplus0})-(\ref{cost0}) is the unique {\it constrained viscosity solution} of the corresponding HJB equation which, even in this simplest case, is a {\it fully nonlinear, non-local}, and {\it degenerate} parabolic partial integro-differential equation.

The main purpose of this paper is to construct the optimal strategy using the solution (whence the value function) of the HJB equation. To describe the main
difficulties in this effort, we begin with the following observation. By simply calculating the maximizer of the Hamiltonian from the HJB equation (see (\ref{HJB0}) below), one can obtain the following candidate of optimal strategy:
\bea
\label{optim}
\left\{\ba{lll}
\g^*_t=\Big[-\frac{(\m-r)V_x(t, X^*_t, W_t)}{\si^2 X^*_t V_{xx}(t, X^*_t, W_t)}\Big]\vee0\wedge 1; \ms\\
a^*_t= \dot L^*_t= M{\bf 1}_{\{V_x(t, X^*_t, W_t)<1\}}+p{\bf 1}_{\{V_{x}(t,X^*_{t},W_{t})=1\}},
\ea\right.
\eea
where $V$ is the viscosity solution and $M>0$ is the given upper bound of the dividend rate, that is, assuming $0\le a_t=\dot L_t\le M$). Then we immediately
see that there are two major technical issues.  First,  the {\it regularity} of the viscosity solution (i.e., the validity of the derivatives $V_x$ and $V_{xx}$), which
 is a  tall order for a non-local, degenerate HJB equation. Second, although the optimal dividend rate does display a ``barrier" nature, the execution time is obviously state-dependent, which raises a serious question about the well-posedness of the resulting closed-loop system. A natural way to get around these
 difficulties is to add some additional Brownian motions to the system so that the corresponding HJB equation becomes non-degenerate, hence possesses classical solutions. Then an argument of ``vanishing viscosity"
might lead to at least some $\e$-optimal strategy. Unfortunately, such a method does not work easily in this model, since the random clock $W$, the key  for the Markovization, cannot be perturbed by a Brownian motion. Therefore the degeneracy of the HJB in the variable $W$ is  unremovable by this approach. Nevertheless, we shall consider an auxiliary HJB-type of PDE, and prove that its solution can be used to
construct the $\e$-optimal strategy without using any control theoretic arguments. Our discussion benefitted greatly from a
recent work on nonlocal HJB equations (cf. \cite{GMS}), except that in the present situation we need to deal with a unbounded domain.

%

The rest of the paper is organized as follows. In section 2  we briefly recall the original problem and introduce all the concepts and notations.
In section 3 we prove the existence and uniqueness of the viscosity solution of our key auxiliary PDE, keeping in mind that such a PDE does not
corresponding to an actual control problem(!).  In Sections 4 we prove the desired convergence of the solutions of the approximating
PDEs to the value function. In Section 5 we construct an prospective $\e$-optimal strategy in terms of the solutions to the approximating PDEs. In Section 6 we prove the well-posedness of the closed-loop system corresponding this strategy, and in Section 7 we verify that the constructed strategy does produce the desired $\e$ optimality.
Some technical results are proved in the Appendix to keep the discussion more readable.

\section{Prelimilaries}

Throughout this paper we consider
a complete probability space $(\O,\cF, \hP)$ on which is
defined  standard Brownian motion $B=\{B_t:t\ge 0\}$, and a {\it renewal}
counting process  $N=\{N_t\}_{t\ge0}$, independent of $B$. More precisely, denoting
$\{\sigma_{n}\}_{n=1}^\infty$ to be the jump times ($\sigma_{0}:=0$) of $N$, and $T_{i}=\sigma_{i}-\sigma_{i-1}$, $i=1,2,\cdots$ to be its waiting times,
we assume that $T_i$'s are independent
and identically distributed,
with a common distribution $F:\hR_+\mapsto \hR_+$. We shall assume that there exists an {\it intensity function} $\l:[0,\infty)\mapsto [0,\infty)$  such that
$\bar F(t):=\hP\{T_1>t\}=\exp\{-\int_0^t\lambda(u)du\}$, so that $\lambda(t)=f(t)/\bar F(t)$, $t\ge 0$, where $f$ is the density
function of $T_i$'s. Clearly, if $\l(t)\equiv \l$ is a constant, then $N$ becomes a standard Poisson process.

Let $T>0$ be a given time horizon, $\hX$ be a generic Euclidean space, and $\cG\subseteq\cF$ be any sub-$\si$-field.
We denote $\hC([0,T];\hX)$ to be the space of continuous functions taking values in $\hX$ with the usual sup-norm;
$L^p(\cG;\hX)$ to be the space of all $\hX$-valued, $\cG$-measurable
random variables $\xi$ such that $\hE|\xi|^p<\infty$,  $1\le p\le\infty$; and
$L^p_\hF([0,T];\hX)$ to be the space of all
$\hX$-valued, $\hF$-progressively measurable processes $\xi$ satisfying
$\hE\int_0^T|\xi_t|^pdt<\infty$, where  $\hF=\{\cF_t:t\ge0\}$ is a given filtration in $\cF$, and $1\le p\le\infty$. Here $p=\infty$ means that
all elements are bounded.

%
%
Given a renewal counting process $N$, we shall consider the following
{\it claim process} for our reserve mode: $Q_t=\sum_{i=1}^{N_t}U_i$, $t\ge 0$, where $\{U_{i}\}_{i=1}^\infty$ is a sequence of
random variables representing the ``size" of the incoming claims. We assume that $\{U_i\}$  are i.i.d. with a common distribution function $G$ 
(and density $g$), 
  independent of $(N, B)$. We note that the process $Q$  is non-Markovian in general (unless the counting process $N$ is a Poisson process), but
  can be ``Markovized" by the so-called {\it Backward Markovization} technique (cf. e.g., 
\cite{RSST}). More precisely,
we define a new process  $W_t=t-\sigma_{N_t}$, $t\ge 0$, that is,
the time elapsed since the last jump. Then  it is known (see, e.g., \cite{RSST}) that
the process $(t, Q_t, W_t)$, $t\ge0$, is a piecewise deterministic Markov process (PDMP). We note that at each jump time $\si_i$, the jump size $|\D W_{\si_i}|=\si_i-\si_{i-1}=T_i$, and $0\le W_t\le t\le T$, $t\in[0,T]$.

Now let us denote  $\{\cF^\xi_t:t\ge 0\}$ to be the natural filtration generated by process $\xi=B, Q, W$, respectively, with the usual $\hP$-augmentation such that it satisfies the {\it usual hypotheses} (cf. e.g.,  \cite{PR}).
Throughout this paper we consider the filtration  $\hF=\hF^{(B,Q, W)}=\{\cF_t\}_{t\ge 0}$, where $\cF_t:=\cF^B_t\vee\cF^Q_t\vee\cF^W_t$, $t\ge 0$.
For any $s\in [0,T]$, let us consider the process $(B,Q, W)$ starting from $s\in[0,T]$.
First assume $W_s=w$, $\hP$-a.s., let us consider the
{\it regular conditional probability distribution} (RCPD) $\hP_{sw}(\cd):=\hP[\,\cd\,|W_s=w]$ on $(\O, \cF)$, and consider the
``shifted" version of processes $(B, Q, W)$ on the space $(\O, \cF, \hP_{sw}; \hF^s)$, where $\hF^s=\{\cF_t\}_{t\ge s}$.
Define $B^s_t:=B_{t}-B_s$, $t\ge s$. Clearly, since $B$ is independent of $(Q, W)$,
$B^s$ is  an $\hF^s$-Brownian motion under $\hP_{sw}$, defined on $[s,T]$, with $B^s_s=0$. Next, we restart the clock at time $s\in[0,T]$ by defining the new counting process $N^s_t:=N_{t}-N_s$, $t\in [s,T]$.  Then, under $\hP_{sw}$, $N^s$ is a ``delayed" renewal process, in the sense that
while its waiting times $T^s_i$, $i\ge 2$, remain i.i.d. as the original $T_i$'s, its ``time-to-first jump", denoted by
$T^{s,w}_1:=T_{N_s+1}-w=\si_{N_s+1}-s$, should have the survival probability
\bea
\label{PTsw}
\hP_{sw}\{T^{s,w}_1>t\}=\hP\{T_1>t+w|T_1>w\}=e^{-\int_w^{w+t}\l(u)du}.
\eea
 In what follows we shall denote $N^s_t\big|_{W_s=w}:=N^{s,w}_t$,
$t\ge s$, to emphasize the dependence on $w$ as well. Correspondingly, we shall denote $Q^{s,w}_t=\sum_{i=1}^{N^{s,w}_t}U_i$
and $W^{s,w}_t:=w+W_{t}-W_s=w+[(t-s)-(\si_{N_{t}}-\si_{N_s})]$, $t\ge s$. It is readily seen that $(B^s_t, Q^{s,w}_t, W^{s,w}_t)$, $t\ge s$, is an $\hF^s$-adapted  process defined on $(\O, \cF, \hP_{sw})$, and it remains Markovian.

\ms
{\bf The Markovized Optimal Investment-Dividend Problem.} Taking the process $W$ into account, we now reformulate the renewal risk model (\ref{surplus0})-(\ref{cost0}) so that
it is Markovian. Similar to our previous work \cite{BaiMa17}, we shall make use of the following {\it Standing Assumptions}:
\begin{assum}
\label{assump0}
(a) The interest rate $r$, the drift $\mu$, the volatility $\si$, and the insurance premium $p$ are all positive constants,;

(b)  The distribution functions $F$ (of $T_i$'s) and $G$ (of $U^i$'s) are continuous on $[0,\infty)$. Furthermore, $F$ is absolutely continuous, with density function $f$ and intensity function $\l(t):=f(t)/\bar F(t)>0$, $t\in[0,T]$;

(c) The cumulative dividend process $L$ is absolutely continuous with respect to the Lebesgue measure. That is,
 there exists $a\in L^2_\hF([0,T]; \hR_+)$, such that $L_t=\int_0^t a_sds$, $t\ge 0$.
 We assume further that for some  constant $M\ge p>0$, it holds that $0\le a_t\le M$, $dt\times d\hP$-a.e.
 \qed
\end{assum}

For any $[s, t]\subseteq [0,T]$, we say that a strategy $\pi=(\g,a )$ is
{\it admissible} on $[s, t]$ if $\pi\in L^2_\hF([s,t];\hR^2)$, such that $\g_u\in[0,1]$, $a_u\in[0,M]$, $u\in[s,t]$, $\hP$-a.s.
More specifically, for any $(s,w)\in [0,T]^2$, we denote the set of all admissible  strategies on $[s,T]$, defined on the probability space $(\O, \cF, \hP_{sw})$
by  $\sU_{ad}^{s,w}[s, T]$. In particular, we denote $\sU^{0,0}_{ad}[0,T]$ by $\sU_{ad}[0,T]=\sU_{ad}$ for simplicity.

Let  $\pi=(\g, a)\in \sU_{ad}^{s,w}[s, T]$, we now consider the ``Markovized" reserve model of (\ref{surplus0})-(\ref{cost0}):
\bea
\label{Xsxw}
\left\{\ba{lll}
\dis dX_{t}=pdt+[r+(\m-r)\g_t]X_{t}dt+\gamma_tX_tdB_t-dQ^{s,w}_t-a_tdt, \q X_s=x;\\
W_t =w+(t-s)-(\si_{N_t}-\si_{N_s}), \qq\qq t\in[s, T],
\ea\right.
\eea
with the expected cumulated dividends up to ruin:
%
\bea
\label{cost1}
J(s,x,w;\pi);=\hE_{sw}\Big\{\int_{s}^{\tau^{\pi}_s\wedge T}e^{-c (t-s)}a_{t}dt\Big|X^\pi_{s}=x\Big\}:=\hE_{sxw}\Big\{\int_{s}^{\tau^{\pi}_s\wedge T}e^{-c (t-s)}a_{t}dt\Big\},
\eea
and the value function:
 \bea
 \label{V1}
 V(s,x,w):= \sup_{\pi\in \sU^{sw}_{ad}[s, T]} J(s,x,w; {\pi}).
 \eea
In the above $(X^\pi,W)=(X^{\pi,s,x,w}, W^{s,w})$ is the solution to (\ref{Xsxw}) and $\t^\pi_s=\t^{\pi,x,w}_s:=\inf\{t >s: X^{\pi,s,x,w}_{t}<0\}$ is the ruin time.

\ms
{\bf The HJB Equation and its Viscosity Solution.} We now briefly recall the main result of \cite{BaiMa17}. We first note that there is a natural domain for the
initial state $(s,x,w)$, denoted by $D:= \{(s,x,w): 0\leq s \leq T, x\geq0, 0\leq w \leq s\}$. Here $w \le s$ is due to the fact $W_t\le t$ always holds. We thus
assume that the value function $V$ is defined on $D$ and that  $V(s,x,w)=0$, for $(s,x,w)\notin D$. We also define the following two sets:
 \bea
 \label{D}
 \sD&:=&\mbox{int}D=\{(s,x,w)\in D: 0<s<T, \, x>0,\, 0<w<s\}; \\
  \sD^*&:=&\{(s,x,w)\in D: 0\leq s<T,\, x\geq0, \, 0\leq w\leq s\}.\nonumber
   \eea
Clearly $\sD\subset \sD^*\subseteq \bar\sD=D$, and $\sD^*$ does not include boundary at the terminal time $s=T$.
Furthermore, we denote $\hC^{1,2,1}_0(D)$ to be the set of all functions $\vf\in \hC^{1,2,1}(\sD)$ such that
for $\eta=\vf$, $\vf_t$, $\vf_x$, $\vf_{xx}$, $\vf_w$, it holds that $\lim_{(t,y,v)\to (s,x,w)\atop  (t,y,v)\in\sD}\eta(t,y,v)=\eta(s,x,w)$,
for all $(s,x,w)\in D$;
and $\vf(s,x,w)=0$, for $(s,x,w)\notin D$.  We note that while
a function $\vf\in\hC^{1,2,1}_0(D)$ is well-defined on $D$, it is not necessarily continuous on
the boundaries $\{(s,x,w): x=0 ~\mbox{or}~ w=0 ~\mbox{or}~ w=s\}$.

Now,  for $\th=(s,x,w)\neg\in\neg D$, $\xi=(\xi^1,\xi^2)\in \hR^2$, $y, A, z\in\hR$, and $(\g,a)\in [0,1]\times [0,M]$, we define the following Hamiltonian:
\bea
\label{H0}
H(\th,y, \xi, A, z, \g, a):=\frac{\sigma^{2}}2\gamma^{2} x^{2}A+[p+(r+(\m-r)\g) x-a]\xi^1
			+\xi^2+ \l(w)z+(a-cy),
\eea
and for $\vf\in\hC^{1,2,1}_0(D)$ we define  the second-order  partial integro-differential operator:
\bea
\label{sL0}
\sL[\vf](s,x,w):=\sup_{\gamma\in [0,1], a\in [0,M]} H(s,x,w, \vf, \nabla \vf, \vf_{xx}, I(\vf), \g, a).
\eea
where  $\nabla\vf:=(\vf_x,\vf_w)$, and  $I[\vf]$ is the integral operator defined by
\bea
\label{Iv}
 I[\vf]:=\int_0^\infty[\vf(s,x-u,0)-\vf(s,x,w)]dG(u)=\int_0^x\vf(s,x-u,0)dG(u)-\vf(s,x,w).
 \eea
Here the last 	equality is due to the fact that $\vf(s,x,w)=0$ for $x<0$.

The main result of \cite{BaiMa17} is that the value function $V$ is the unique {\it constrained} viscosity solution of the following
HJB equation:
		\bea
			\label{HJB0}
			\{V_s+\sL[V]\}(s,x,w)=0; \q (s,x,w)\in \sD;\qq
			V(T, x, w)=0.
			\eea
To facilitate our future discussion, we end this section by recalling the definition of the ``constrained viscosity solution" to the {\it Partial Integro-Differential Equation} (PIDE) (\ref{HJB0}) (cf. \cite{BaiMa17}):
 \begin{defn}
 \label{defvis}
		Let ${\cal O}\subseteq \sD^*$ be a subset such that  $\pa_T\cO:=\{(T,y,v)\in\pa\cO\}\neq \emptyset$, and let $v\in \hC(\cO)$.
		
		(a) We say that $v$ is a viscosity subsolution of (\ref{HJB0}) on ${\cal O}$,   if $v(T,y,v)\le 0$, for $(T, y, v)\in \pa_T\cO$;
 and for any $(s,x,w)\in \cO$ and $\vf\in \hC^{1,2,1}_0(\cO)$
such that $0=[v-\vf](s,x,w)=\max_{(t,y,v)\in\cO}[v-\vf](t,y,v)$, it holds that
\bea
\label{vsub}	
			\vf_s(s,x,w)+\sL[\vf](s,x,w)\geq 0.
\eea

		(b) We say that $v$ is a viscosity supersolution of (\ref{HJB0}) on ${\cal O}$,   if $v(T,y,v)\ge 0$, for $(T, y, v)\in \pa_T\cO$;
 and for any $(s,x,w)\in\cO  $ and $\vf\in \hC^{1,2,1}_0(\cO)$ such that $0=[v-\p](s,x,w)=\min_{(t,y,v)\in\cO}[v-\vf](t,y,v)$,
it holds that
\bea
\label{vsup}	
			\vf_s(s,x,w)+\sL[\vf](s,x,w)\leq 0.
\eea
		
		(c) We say that $v\in \hC(D)$ is a ``constrained viscosity solution" of (\ref{HJB0}) on $\sD^*$ if it is both a  viscosity subsolution of (\ref{HJB0}) on $\sD^*$ and a viscosity supersolution of (\ref{HJB0}) on ${\sD}$.
\qed
	\end{defn}

\section{An Auxiliary  Equation}
\setcounter{equation}{0}

As we pointed out, our goal is to construct a sensible approximation of the optimal strategy based on the explicit form (\ref{optim}) using the solution to the HJB equation (\ref{HJB0}). But the degenerate nature of the Hamiltonian (\ref{H0}), especially in the variable $w$, makes this task particularly challenging, since the random clock $W=\{W_t\}$ cannot be perturbed by a Brownian motion, in order to keep the Markovization procedure intact. As a remedy in the
rest of this paper we shall therefore focus mostly on the PDE aspect of the issue and introduce an auxiliary non-degenerate PDE that is of the same structure as (\ref{HJB0}), but its solution cannot be regarded as a value function to any stochastic control problem. As a consequence our arguments could be more
analytical than some of the control theoretical ones in the literature, but they are interesting in its own right. In fact, to the best of our knowledge, the
regularity of the constrained viscosity solution to a non-local HJB equation of this particular type on a unbounded domain is new.

Our plan of attack is quite similar to that of the recent work \cite{MS}. More precisely, we begin with the following extended domain of $D$: for each $\d>0$,
\bea
\label{Dd}
D_{\delta}=\{(s,x,w): 0< s\leq T+\delta, x\geq-\delta, -\delta\leq w \leq s+\delta\}.
\eea
As before, we denote $\sD_\d:=$ int$D_\d$, and consider the  ``truncated" complement" of $D_\d$:
\bea
\label{Ddc}
\sD^{*,c}_{\delta}:=(\{T+\delta\}\times \hR^{2})\cup \big(\cup_{0<s< T+\delta} \sD^c_{\d, s}\big),
\eea
where for $0<s<T+\d$, $\sD_{\delta,s}=\{(x,w): x>-\delta, -\delta < w < s+\delta\}$ is the $s$-section of $\sD_\d$, and $\sD^c_{\d, s}$ is the complement
of $\sD_{\d, s}$.
Clearly, $D_\d\cup \sD^{*,c}_\d=(0,T+\d]\times\hR^2$.

Next, we define a ``perturbed" non-degenerate Hamiltonian.  Let $\e_n>0$, $n=1,2,\cds$ be a sequence such
that $\e_n\da 0$, as $n\to\infty$. We define for  $\th=(s,x,w)\in D_\d$, $\xi=(\xi^1,\xi^2)\in \hR^2$, $y,  A_1, A_2, z\in\hR$, and $(\g,a)\in [0,1]\times [0,M]$,
\bea
\label{Hn}
H^{n}(\th,y,\xi,A_1,A_2,z, \g, a):= H(\th, y, \xi, A_1, z, \g, a)+\frac{\varepsilon_{n}}{2}A_1+\frac{\varepsilon_{n}}{2}A_2,
\eea
where $H$ is the Hamiltonian defined by (\ref{H0}), and consider the following auxiliary PIDE:
 \bea
\label{HJB2}
\left\{\ba{lll}
v_{t}(s,x,w)+\sL^{n,\d}[v](s,x,w)=0, \qq
\mbox{on  $ \sD_\d$,}\\
v(s,x,w)=\Psi(s,x,w),\qq (s,x,w)\in \sD^{*,c}_{\delta}.
\ea\right.
\eea
Here, as before, for a smooth function $\vf$ and $\nabla \vf=(\vf_x, \vf_w)$,
\bea
\label{Lnd}
\left\{\ba{lll}
\dis \sL^{n,\d}[\vf](s,x,w):=\sup_{\gamma\in [0,1], a\in [0,M]} H^{n}(s,x,w,\vf,\nabla \vf, \vf_{xx}, \vf_{ww},I^{\delta}[\vf],\g, a),  \ms\\
 \dis I^{\delta}[\vf](s,x,w):=\int_0^{x+\delta}\vf(s,x-u,-\delta)dG(u)-\vf(s,x,w);
 \ea\right.
 \eea
and $\Psi$ is a function to be determined later.  
We shall argue that there exists a unique classical solution to (\ref{HJB2}), denoted  by $V^{n,\d}$,  such that $\lim_{n\to\infty, \d\to 0}V^{n,\d}=V$,  the value function $V$ defined by (\ref{V1}), uniformly on compacta.

We should note that since the equation (\ref{HJB2}) does not necessarily correspond to any stochastic control problem, the existence of the solution, even in the
viscosity sense, is not clear. In the rest of this section we shall first show that there is indeed a viscosity solution to this equation, and in the next section
we shall argue that such a solution is actually the unique classical solution. To simplify the argument we shall assume $0<\delta<1$ throughout our discussion.

{\bf The function $\Psi$.} We now give a detailed description of the function $\Psi$, which is crucial for our construction of the viscosity solution. We first note
that once such a function is chosen, we can modify the PIDE (\ref{HJB2})  to one with homogeneous boundary condition via the following  standard transformation. Assume that $\Psi$ is a (smooth) boundary condition. Let $\tilde v=v-\Psi$, then we have
\bea
\label{HJB2E}
\left\{\ba{lll}
 (\tilde v+\Psi)_{t}+\sL^{n,\d}[\tilde v+\Psi]= v_t+\sL^{n,\d}_\Psi[\tilde v]=0,\\
\tilde v(s,x,w)=0,\ \ (s,x,w)\in \sD^{*,c}_{\delta}.
\ea\right.
\eea
where
$\sL^{n,\d}_{\Psi}[\vf]:=\Psi_t+ \sL^{n,\d}[\vf+\Psi]$ will have the  same properties as $\sL^{n, \d}$.
Furthermore, we shall make the following assumption. Recall the set $D_\d$ and the constants $M>0$ in Assumption \ref{assump0}.
\begin{assum}
\label{assump1} There exists $\Psi\in \hC^{1,3,3}(\hR^3)$ such that

 (i)  there exists $K_1>0$ such that $0\le \Psi(\th)\le K_1$, $\th=(s,x,w) \in D_{1}$, and $\Psi(\th)=0 $,  $\th\in D^{c}_{1}$;

(ii)  there exists  $0<K_{2}<M$, such that for  any  $\th\in D_{1}$,
\beaa
M-K_{2}\leq \Psi_{t}+H^{n}(\th, \Psi,\nabla \Psi, \Psi_{xx},\Psi_{ww},I^{\delta}[\Psi],0, M), \q 0<\delta<1,\q n\geq 1;
\eeaa

(iii)  $\Psi(s,x,w)$ is strictly increasing with respect to $x$, and  for some $0<\d_0<1$,
 \bea
 \label{b}
b :=\inf_{(s,x,w)\in (0,T]\times  [-\d_0,0]\times [0,s]}\Psi_{x}(s,x,w)>1.
\eea
 \end{assum}
%
%
%

We should note that under Assumption \ref{assump0}, Assumption \ref{assump1}-(ii) holds if $M$ is large enough,  but
(iii) is a special requirement that is important in our convergence analysis.
In the rest of the paper we shall fix a function $\P$ satisfying Assumption \ref{assump1}, and consider viscosity solution within a special class
of functions associated to $\P$.
More precisely, we have the following definition.
\begin{defn}
\label{classP}
We say that a function $v$ is of class \mbox{\rm ($\P$)} if it satisfies the following conditions:

\ss
(1)  $v(s,x,w)=\Psi(s,x,w),\ \ (s,x,w)\in \sD^{*,c}_{\delta}$;

(2) $v(s,x,w)$ is  increasing with respect to $x$ on $D_{\delta}$;

(3)  $v(s,x,w)$ is  bounded on $ D_{\delta}$;

(4)  $v(s,x,w)-v(s,-\delta,w)\geq x+\delta$ as $x\downarrow -\delta$ for any $0\leq s \leq T+\delta, -\delta\leq w \leq s+\delta$.
\qed
\end{defn}

We shall construct a viscosity solution  of (\ref{HJB2})
that is of class ($\P$) by the well-known Perron's method. To begin with, we need an
important lemma, whose proof will be deferred to the Appendix  in order not to disturb the flow of discussion. \begin{lem}
\label{classL1}
Assume Assumptions \ref{assump0} and \ref{assump1}. There exist both viscosity supersolution $\overline{\psi}$ and subsolution $\underline{\psi}$
of class \mbox{\rm ($\P$)} to (\ref{HJB2}) on $\sD_{\delta}$. Furthermore, it holds that $\overline{\psi}=\underline{\psi}=\Psi$ on $\sD^{*, c}_{\delta}$.
\qed
\end{lem}

Next, for given $\Psi$, we consider the following set
$${\sF}= \{v: v \ \mbox{is a  viscosity subsolution of class {($\P$)} to (\ref{HJB2}) on $\sD_{\delta}$}, \mbox{ s.t. } \underline{\psi}\leq v \leq \bar{\psi} \},$$
where $\underline{\psi}$ and $\bar{\psi}$ are  the viscosity subsolution and supersolution, respectively, of class ($\P$) mentioned in Lemma
\ref{classL1}.
Define
\bea
\label{u0}
u(s,x,w):= \sup_{ v\in {\sF}} v(s,x,w), \qq (s,x,w)\in D_\d,
 \eea
and let $u^*$ (resp. $u_*$) be the {\it upper semicontinous (USC) envelope} (resp. {\it lower semicontinous (LSC) envelope}) of $u$, defined
respectively by
 \bea
 \label{u1}
\left\{\ba{lll}
\dis u^{*}(s,x,w):=\lim_{r\da 0} \sup\{u(t,y,v): (t,y,v)\in \sD_{\delta}, \   \sqrt{|t-s|^{2}+|y-x|^{2}+|v-w|^{2}} \leq r  \}, \ms\\
\dis u_{*}(s,x,w):=\lim_{r\da 0} \inf\{u(t,y,v): (t,y,v)\in \sD_{\delta}, \   \sqrt{|t-s|^{2}+|y-x|^{2}+|v-w|^{2}} \leq r  \}.
\ea\right.
\eea
 The main result of this section is the following theorem, which obviously implies the existence of the viscosity solution to (\ref{HJB2}).

\begin{thm}
\label{060501}
Assume that Assumptions \ref{assump0} and \ref{assump1} are in force. Then $u^{*}$ (resp. $u_{*}$) is  a viscosity subsolution
(resp. supersolution) of class \mbox{\rm ($\P$)} to (\ref{HJB2}) on $\sD_{\delta}$.
\end{thm}

{\it Proof.} The fact that $u^*$ is a subsolution is more or less straightforward, we shall omit the proof and accept it as a fact, and
 prove only  that $u_*$ is a supersolution of class ($\P$).  It is easy to verify that $u_{*}$   of class ($\P$).
Suppose that $u_{*}$ is not a supersolution, then there exists $\th_0=(s_{0},x_{0},w_{0})\in \sD_{\delta}$ and $\varphi\in \hC^{1,2,2}_0(\sD_\d)$ such that $0=[u_{*}-\varphi](\th_0)<[u_*-\vf](\th)$, for all $\th\in\sD_\d$, but
\beaa
\label{notsupsol0}
\pa_{t}\varphi(\th_{0})+\sup_{\gamma\in [0,1], a\in [0,M]}H^{n}(\th_{0}, u_{*},\nabla\vf,  \varphi_{xx}, \varphi_{ww},I^{\delta}[\vf], \g, a)=: \varepsilon_{0}>0.
\eeaa
By continuity, we can then find $\eta_{0}>0$ such that, for any $\th\in B_{\eta_0}(\th_{0})\subset \sD_{\delta}$,
\bea
\label{notsupsol}
\pa_{t}\varphi(\th)+\sup_{\gamma\in [0,1], a\in [0,M]}H^{n}(\th, u_{*},\nabla\vf,  \varphi_{xx}, \varphi_{ww},I^{\delta}[\vf], \g, a)> \varepsilon_{0}/4.
\eea

We shall argue that (\ref{notsupsol}) means that one can construct a subsolution $\psi^*\in \sF$, such that $\psi^*(\th_0)>u(\th_0)$, which would contradict the definition of $u$. To this end,  note that being of  class ($\P$) $u_*$ is increasing in $x$. Thus for $0<\e_{1}<\frac{\e_{0}}{2}$, we can modify $\vf$ slightly so that
on $B_{\eta_{0}}(\th_0)$ (or choose a smaller ball if necessary) $\vf$ is increasing in $x$,  but it is
decreasing in $x$ for $x$ sufficiently large, such that
 \bea
\label{101502}
\dis \inf_{\th\in B^{c}_{\eta_{0}}(\th_{0})\cap D_{\delta}} \{u_{*}(\th)-\vf(\th)\}\geq \varepsilon_{1}>0.
\eea
Note that, by definition of  $u$, we have  $\vf\leq u_{*} \leq \bar{\psi}$ in $\sD_{\delta}$. We claim that $\vf(\th_{0})<\bar{\psi}(\th_{0})$. Indeed, if $\vf(\th_{0})=u_{*}(\th_{0})=\bar{\psi}(\th_{0})$, then $\bar{\psi}-\vf$ has a strict minimum at $\th_{0}$. Since $\bar{\psi}$
is a viscosity supersolution (\ref{HJB2}) on $\sD_{\delta}$, we have
$$\pa_t\vf(\th_{0})+\sup_{\g\in [0,1], a\in [0,M]} H^n(\th_{0}, \vf,\nabla\vf, \pa_{xx}\vf,\pa_{ww}\vf,I^\d[\vf],
\g, a)\leq 0,$$
contradicting  (\ref{notsupsol}). Therefore, by continuity of $\bar{\psi}$ and $\vf$, we can find  $0<\eta_{2}<\eta_{0}$ and
$\e_{2}>0$, such that  $\vf(\th)<\bar{\psi}(\th)-\e_{2}$,  $\th\in B_{\eta_{2}}(\th_0)$. Note that   $u_{*}-\vf$ has a strict minimum at
$\th_0$,  we have
\bea
\label{Deltar}
\Delta_{r}:=\inf_{\th\in B_{r}^{c}(\th_{0})\cap D_{\delta}}\{u_{*}(\th)-\vf(\th)\}=\inf_{\th\in \ol{B_{r}^{c}}(\th_{0})\cap \ol{D}_{\delta}}\{u_{*}(\th)-\vf(\th)\}>0, \qq r>0.
\eea

Let us now fix $r_0\in (0, \eta_2)$. Recall that we have modified $\vf$ so that for some $\hat x>0$ large enough, it is  decreasing in $x$, for $x>\hat x$. We assume without loss of
generality that $\hat x>x_0+r_0$. Define
$E_{\delta}(\hat x):=\{\hat\th:=(s,\hat x, w): 0\leq s<T+\delta, -\delta< w <s+\delta\}$.
Clearly, $E_\d(\hat x)\subset \ol{B^c}_{r_0}\cap \ol{D}_\d$, thus by (\ref{Deltar}) we have
$u_{*}(\hat \th)-\vf(\hat\th)\geq \Delta_{r_0}$, for  $\hat\th \in E_\d(\hat x)$.
Now for fixed $\hat \th_1=(s_1, \hat x, w_1)\in E_\d(\hat x)$, by  definition of $u_{*}$ we can choose  $\hat v_1\in {\sF}$ such that
$\hat v_1(\hat\th_1)-\vf(\hat\th_1)\geq \frac{3\Delta_{r_0}}{4}$.
But since $\hat v_1\in\sF$ (whence increasing in $x$) and $\vf $ is decreasing in $x$ for $x>\hat{x}$, we have
\bea
\label{hatv-vf}
\hat v_1(s_1,x,w_1)-\vf(s_1,x,w_1)\geq \hat v_1(\hat\th_1)-\vf(\hat\th_1)\ge\frac{3\Delta_{r_0}}{4}, \qq \mbox{for $x\geq \hat x $.}
\eea
On the other hand, by continuity of $(\hat v_1-\vf)(\cdot,\hat x,\cdot)$, there exists $\hat\eta_1>0$, such that
\bea
\label{hatv-vf1}
\inf_{(s,w)\in \bar{B}_{\hat \eta_1}(s_1, w_1)\cap \bar{E}_{\delta}(\hat x)}\{\hat v_1(s,\hat x,w)-\vf(s,\hat x,w)\}\geq \frac{\Delta_{r_0}}{2}.
\eea
Note that $\bar{E_{\delta}}(\hat x)$ is compact,  there exists a finite set $\{(s_{j},w_{j})\}_{j=1}^{m_0}\subset \bar{E_{\delta}}(\hat x)$, together
with $\hat v_j\in\sF$  and constants
$\hat\eta_j>0$, $j=1, \cds, m_0$, such that $\bar{E_{\delta}}(\hat x)\subset \cup_{j=1}^{m_{r}}\bar{B}_{\hat\eta_j}(s_{j},w_{j})$, and both (\ref{hatv-vf}) and
(\ref{hatv-vf1}) hold for each $j$. Now let us define
$$\ell_{0}(\th)=\sup_{1\leq j \leq m_0}\hat v_j(\th), \qq\q   \th \in \sD_{\delta}.
$$
Then one can check, as before, that  $\ell_0\in {\cal F}$, and is increasing with $x$ on $\sD_{\delta}$. Furthermore, since each
$\hat v_j$ satisfies (\ref{hatv-vf}) and (\ref{hatv-vf1}), it is readily seen that
\bea
\label{l0-vf}
\inf_{(s,x,w)\in \sD_{\delta}\backslash D_{\delta,\hat x}  }\{\ell_0(s,x,w)-\vf(s,x,w)\}\geq \frac{\Delta_{r_0}}{2},
\eea
where $D_{\delta,\hat x}:=\{(s,x,w): 0< s<T+\delta, -\delta<x<\hat x, -\delta< w <s+\delta\}$.

Now let us consider the set $\bar {D}_{\delta,\hat x}\backslash B_{r_0}(\th_0)$. By (\ref{Deltar}) we have
$u_{*}(\th)-\vf(\th)\geq \Delta_{r_0}$ for all $\th\in \bar {D}_{\delta,\hat x}\setminus B_{r_0}(\th_0)$. Since $\bar {D}_{\delta,\hat x}\backslash B_{r_0}(\th_0)$
is compact, we can repeat the same argument as before
to obtain a $\ell_1\in\sF$ so that
\bea
\label{l1-vf}
\inf_{(s,x,w)\in \overline{D}_{\delta,\hat x}\setminus B_{r_0}(\th_{0}) }\{\ell_1(s,x,w)-\vf(s,x,w)\}\geq \frac{\Delta_{r_0}}{2}.
\eea

Let $0<\a_0<\min\{\frac{\e_2}{\Delta_{r_0}}, \frac12\}$, and define
\bea
U(\th):=\left\{\ba{lll}
\max\{\vf(\th)+\a_0 \Delta_{r_0}, \ell_0(\th),\ell_1(\th) \},  \qq \qq & \mbox{if $ \th\in  B_{r_0}(\th_{0})$}\\
\max\{\ell_0(\th),\ell_1(\th)\} & \mbox{if  $\th\in B_{r_0}^{c}(\th_{0})\cap \sD_\d$,}
\ea\right.
\eea
Then, by the choice of $r_0$ and $\a_0$, we have $\underline{\psi}\leq U \leq \bar{\psi}$ in
 $\sD_{\delta} $,  and
 \bea
 \label{Ugeu}
U(\th_{0})\geq \vf(\th_{0})+\a_0 \Delta_{r_0}>\vf(\th_{0})=u_{*} (\th_{0}).
\eea
%
%
We claim that $U$ is a viscosity subsolution of class ($\P$) to (\ref{HJB2}) in $\sD_{\delta}$, which would be a contradiction to the the definition
of $u_*$ and prove the theorem.

To this end, For any $\bar\th:=(t,y,v)\in \sD_{\delta}$, suppose that there is a function $\phi\in  \hC^{1,2,2}_0(\sD_{\delta})$ such that $0=U(\bar\th)-\phi(\bar\th)$ is a strict maximum over $\sD_{\delta}$. Consider two possible cases:

{\it Case 1}: $U(\bar\th)=\ell_0(\bar\th)$ or $\ell_1(\bar\th)$. ~
We shall only consider the case $U(\bar\th)=\ell_0(\bar\th)$, as the other case is similar.   Since $\ell_0\leq U \leq \phi$ on $\sD_{\delta}$,  $\ell_0-\phi$ has a maximum  at $\bar\th$. Recall again that, as the ``sup" of subsolutions, $\ell_0$ is
a viscosity subsolution  of (\ref{HJB2}) on $\sD_{\delta}$ as well,  hence we  have
\bea
\label{phieq}
\pa_t\phi(\bar\th)+\sup_{\gamma\in [0,1], a\in [0,M]} H^{n}(\bar\th, \phi,\nabla\phi, \phi_{xx},\phi_{ww},I^{\delta}[\phi], \g, a)\geq 0.
\eea

{\it Case 2}: $U(\bar\th)=\vf(\bar\th)+\a_0 \Delta_{r_0}$. In this case we must have $\bar\th\in B_{r_0}(\th_0)$ by definition of $U$. But since $\vf+\a_0\Delta_{r_0}\leq U \leq \phi$ in $B_{r_0}(\th_{0})$ by our choices of $r_0$ and $\a_0$,  we have $\vf+\a_0\Delta_{r_0} - \phi\leq 0$ in  $B_{r_0}(\th_{0})$. On the other hand, note that
  $\phi\geq U=\max\{\ell_0,\ell_1\}$ in $B_{r_0}^{c}(\th_{0})\cap \sD_{\delta}$, we conclude that
  $$\vf+\a_0 \Delta_{r_0}- \phi\leq \vf+\a_0 \Delta_{r_0}- \max\{\ell_0,\ell_1\}\leq -\frac{\Delta_{r_0}}{2}+\a_0 \Delta_{r_0}\leq 0, $$
in $B_{r_0}^{c}(\th_{0})\cap \sD_{\delta}$. That is,  $\vf+\a_0\Delta_{r_0}-\phi$ has a maximum  at
$\bar\th\in B_{r_0}(\th_{0})\subset B_{\eta_{1}}(\th_{0})$. Then, by (\ref{notsupsol}), choosing $\a_0$ sufficiently small if necessary we have
\bea
\label{phieq1}
&&\pa_t\phi(\bar\th)+\sup_{\gamma\in [0,1], a\in [0,M]}H^{n}(\bar\th,\phi,\nabla\phi, \phi_{xx},\phi_{ww},I^{\delta}[\phi], \g, a)\\
&&\geq \pa_t\vf(\bar\th)+\sup_{\gamma\in [0,1], a\in [0,M]}H^{n}(\bar\th, \vf+\a_0 \Delta_{r_0},\nabla\vf, \vf_{xx},\vf_{ww},I^{\delta}[\vf+\a_0 \Delta_{r_0}], \g, a)\geq 0. \nonumber
\eea
Combining (\ref{phieq}) and (\ref{phieq1}) we conclude that $U$ is  a viscosity subsolution of class ($\P$) to (\ref{HJB2}) in $\sD_{\delta}$, and
$U(\th_0)>u(\th_0)$, a contradiction. This proves the theorem.
\qed
%
%
%
%

Let us now denote the solution to (\ref{HJB2}) by $V^{n,\d}$. We shall argue that such a viscosity solution is unique, and is actually a classical
solution. The proof of uniqueness  will depend on the {\it comparison theorem} as usual, and in this case it can be argued along
the same lines of that in \cite{BaiMa17}, except for some slight modifications. We shall give only a sketch of the proof for completeness.
\begin{thm}[Comparison Principle]
\label{comp2}
Let  $\bar{u}$ be a viscosity supersolution and $\underline{u}$ be a viscosity subsolution of (\ref{HJB2}) on $\sD_{\delta}$, and both
are of class \neg {\rm ($\P$)}. Then
$\underline{u}\leq \bar{u}$ on $D_{\delta}$.
	Consequently,
$u^{*}=u_{*}=:u$ defined by (\ref{u1}) is a unique continuous viscosity solution of class \neg{\rm ($\P$)} to (\ref{HJB2}).
\end{thm}

{\it Proof.}  We first perturb the supersolution slightly so that all the inequalities involved become strict. To this end we define, for $\rho>1, \vartheta>0$, $\bar{u}^{\rho}(t,y, v)=\rho\bar{u}(t,y, v)+\frac{\vartheta}{t}$. Then it is straightforward to check that $\bar{u}^{\rho,\vartheta}(t,y, v)$ is also a supersolution of (\ref{HJB2}) on ${\sD}_{\delta}$ (see, e.g., \cite{BaiMa17}). Furthermore, it is readily seen that $\lim_{t\rightarrow 0}	\bar{u}^{\rho,\vartheta}(t,y,v)=+\infty$; and
\bea
\label{uleP}
\underline{u}(t, y, v)=\Psi(t,y,v)< \rho\Psi(t,y,v)+\frac{\vartheta}{t}= \bar{u}^{\rho,\vartheta}(t,y,v), \qq (t,y,v)\in  \sD^{*,c}_{\delta}.
\eea

   We shall argue that $\underline{u}\leq\bar{u}^{\rho,\vartheta}$ on $\sD_{\delta}$, which, together with (\ref{uleP}), would imply that
$\underline{u}\leq\bar{u}^{\rho,\vartheta}$ on $D_{\delta}$, hence the desired
comparison result as $\lim_{\rho\da 1,\vartheta \da 0}\bar{u}^{\rho,\vartheta}=\bar u$.

We prove this by contradiction. Suppose not, then there exists $\th_0=(t_{0}, y_{0}, v_{0})\in \sD_{\delta}$ such that
$\underline{u}(\th_{0})-\bar{u}^{\rho,\vartheta}(\th_{0})=2\vartheta_{1}>0$, for some $\vartheta_1>0$. Let $\e>0$ be such that
\bea
\phi(\th_{0}):=\underline{u}(\th_{0})-\bar{u}^{\rho,\vartheta}(\th_{0})-2\e (y_{0}+v_{0}+2\d)\geq\vartheta_{1}>0,
\eea
and $\th_\e=(t_{\e}, y_{\e}, v_{\e})\in \bar{{\sD}_{\delta}}$ be such that
\beaa
	\label{Mb}
	M_\e&:=&\sup_{\th=(t, y,v)\in\sD_\delta}(\underline{u}(\th)-\bar{u}^{\rho,\vartheta}(\th)-2\e (y+v+2\d))=\underline{u}(\th_{\epsilon})-
	\bar{u}^{\rho,\vartheta}(\th_{\e})-2\e(y_{\e}+w_{\e}+2\d)\\
&\geq& \underline{u}(\th_{0})-\bar{u}^{\rho,\vartheta}(\th_{0})-2\e (y_{0}+v_{0}+2\d)\geq \vartheta_{1}>0.
		\eeaa
Since $\underline{u}\leq \bar{u}^{\rho,\vartheta}$ on $\partial D_{\delta}$ we see that $\th_{\e}\in \sD_{\delta} $.
Next, for  $\e>0$, we define an auxiliary function: for $\Theta:=(t,x,w,y,v)\in\sC_0:=\{\Th=(t,x,w,y,v): t\in[0,T+\delta], x,y\in[-\delta,+\infty), w,v\in[-\delta,t+\delta]\}$,
	\bea
	\label{Sigma0}
		 \Sigma_{\varsigma,\e}(\Th)=\underline{u}(t,x,w)-\bar{u}^{\rho,\vartheta}(t,y,v)-\e(x+w+y+v+4\delta)-\frac{1}{2\varsigma}(x-y)^{2}-\frac{1}{2
		 \varsigma} (w-v)^{2}.
	\eea
Now let us fix $\e>0$. Let $\Th_\vsi=(t_\vsi, x_\vsi,w_\vsi,y_\vsi,v_\vsi)\in \sC_0$ be such that
\beaa
\label{Me}
M_{\vsi,\e}:=\max_{\Th\in \sC_0}\Sigma_{\vsi,\e}(\Th)=\Sigma_{\vsi,\e}(\Th_\vsi).
\eeaa
By a standard argument (cf. e.g,  \cite{CIL} or \cite{BaiMa17}), using the fact that  $\underline{u}$ is USC and bounded on $D$, and that
$\Sigma_{\vsi,\e}(\Th_\vsi)= M_{\vsi,\e}\geq \Sigma_{\vsi,\e}(t_{\e}, y_{\e}, v_{\e},  y_{\e}, v_{\e})=M_{\e}>0$,
it is not hard to show that there exists $\vsi_0>0$, such that
$\Th_\vsi\in$ int$\,\sC_0$, whenever $0<\vsi<\vsi_0$.

Now applying \cite[Theorem 8.3]{CIL} and following some standard arguments using the  equivalent definition of viscosity solutions in terms of the ``super-jets" (see \cite{CIL}), 
one shows that for any $\delta_{2}>0$, there exist $q=\hat{q}\in \hR$ and  symmetric matrices $A=[A_{ij}]_{i,j=1}^2$ and $B=[B_{ij}]_{i,j=1}^2$,  such that
\beaa
		A_{11}x_{\vsi}^{2}-B_{11}y_{\vsi}^{2}\leq \frac{3}{\vsi}(x_{\vsi}-y_{\vsi})^{2},
	\eeaa
and with $\th_\vsi:=(t_{\vsi},x_{\vsi},w_{\vsi})$, $\bar\th_\vsi:=(t_{\vsi},y_{\vsi},v_{\vsi})$, $\xi^{1, \e}_{\vsi}:=( (x_{\vsi}-y_{\vsi})/\vsi+\e, (w_{\vsi}-v_{\vsi})/\vsi+\e)$ and
 $\xi^{2, \e}_{\vsi}:=( (x_{\vsi}-y_{\vsi})/\vsi-{ \e}, (w_{\vsi}-v_{\vsi})/\vsi-{ \e})$, it holds that
	$$
\left\{\ba{lll}
\dis q+\sup_{\gamma\in [0,1], a\in [0,M]}H^{n}(\th_\vsi, \underline{u},\xi^{1,\e}_{\vsi},A_{11}, A_{22}, I^{\delta}[\underline{u}], \g, a) \geq 0
\\
\dis q+\sup_{\gamma\in [0,1], a\in [0,M]}H^{n}(\bar\th_\vsi,\bar{u}^{\rho,\vartheta},\xi^{2,\e}_{\vsi}, B_{11}, B_{22}, I^{\delta}[\bar{u}^{\rho,\vartheta}], \g, a) \leq 0.
\ea\right.
	$$
	Thus, if we choose $(\gamma_{\vsi},a_{\vsi})\in {\rm arg max}_{(\gamma,a)\in [0,1]\times [0,M]}H^{n}(\th_\vsi,\underline{u},\xi^{1,\e}_{\vsi},A_{11}, A_{22}, I^{\delta}[\underline{u}], \g, a)$, then we have
	\begin{eqnarray*}
	H^{n}(\th_\vsi,\underline{u},\xi^{1,\e}_{\vsi},A_{11}, A_{22}, I^{\delta}[\underline{u}], \gamma_{\vsi},a_{\vsi})
	 -H^{n}(\bar\th_\vsi,\bar{u}^{\rho,\vartheta},\xi^{2,\e}_{\vsi}, B_{11}, B_{22}, I^{\delta}[\bar{u}^{\rho,\vartheta}], \gamma_{\vsi},a_{\vsi})\geq 0.
	\end{eqnarray*}
In other words, by definitions (\ref{H0}) and (\ref{Hn}), this amounts to saying that
	\begin{eqnarray}
	\label{91802}
		&&c(\underline{u}(\th_\vsi)-\bar{u}^{\rho,\vartheta}(\bar\th_\vsi))+\l(w_{\vsi})\underline{u}(\th_\vsi)-
		\l(v_{\vsi})\bar{u}^{\rho,\vth}(\bar\th_{\vsi})\nonumber\\
		&\leq &\frac{1}{2}\sigma^{2}{\g_{\vsi}}^{2}(A_{11}x_{\vsi}^{2}-B_{11}y_{\vsi}^{2})+\frac12\sum_{i=1}^2 (A_{ii}-B_{ii})\e_n+[r+(\m-r)\g_\vsi]\frac{(x_{\vsi}-y_{\vsi})^{2}}{\vsi}+ 2\e( p-a_{\vsi})\nonumber \\
		&&+2\e+\l(w_{\vsi})\int_0^{x_{\vsi}+\delta} \underline{u}(t_{\vsi},x_{\vsi}-u,-\delta)dG(u)-\l(v_{\vsi})\int_{0}^{y_{\vsi}+\delta} \bar{u}^{\rho,\vartheta}(t_{\vsi},y_{\vsi}-u,-\delta)dG(u) \\
		&\leq &\Big(\frac{3\sigma^{2}}2 	+\m\Big)\frac{(x_{\vsi}-y_{\vsi})^{2}}{\e}+\frac12\sum_{i=1}^2 (A_{ii}-B_{ii})\e_n+{ 2\e ( p-a_{\vsi})+2\e} \nonumber \\
		&&+\l(w_{\vsi})\int_{0}^{x_{\vsi}+\delta} \underline{u}(t_{\vsi},x_{\vsi}-u,-\delta)dG(u)-\l(v_{\vsi})\int_{0}^{y_{\vsi}+\delta} \bar{u}^{\rho,\vartheta}(t_{\vsi},y_{\vsi}-u,-\delta)dG(u).\nonumber
	\end{eqnarray}
Now, we can find a sequence $\vsi_{m}\rightarrow 0$ such that
	$\Th_{\vsi_m}:=(t_{\vsi_{m}}, x_{\vsi_{m}},w_{\vsi_{m}}, y_{\vsi_{m}},v_{\vsi_{m}})\rightarrow (\hat{t},\hat{x},\hat{w},\hat{y},\hat{v})\in \bar{\sC_0}$ (here we allow $\hat x=\infty$). Then, a similar argument as before one shows that $(\hat{t},\hat{x},\hat{w},\hat{y},\hat{v}) \notin \pa\sC_0$, $\hat{w}=\hat{v}$,  $\hat x=\hat y<+\infty$, and
\bea
\label{Mep}
		\underline{u}(\hat{t},\hat{x},\hat{w})-\bar{u}^{\rho,\vth}(\hat{t}, \hat{x},\hat{w})-2\e(\hat x+\hat w+2\delta)\geq\lim_{m\to\infty}
\Sigma_{\vsi_{m},\e}(t_{\e}, y_{\e}, v_{\e};  y_{\e}, v_{\e})\geq M_\e>0,
\eea
Thus $\hat\th:=(\hat t,\hat x,\hat w)\in \sD_\delta$. Replacing $\vsi$ by $\vsi_{m}$ and letting ${m}\to\infty$ in
(\ref{91802}), we see from (\ref{Mep}) that
	\begin{eqnarray*}
&&	(c+\l(\hat{w}))(M_{\e}+2\e(\hat x+\hat w+2 \delta))	\le (c+\l(\hat{w}))(\underline{u}(\hat\th)-\bar{u}^{\rho,\vth}(\hat\th))\\
\neg&\neg\leq\neg &\neg\neg \l(\hat{w})\int_{0}^{\hat{x}+\delta}[\underline{u}(\hat{t},\hat{x}-u,-\delta)- \bar{u}^{\rho,\vth}(\hat{t},\hat{x}-u,-\delta)]dG(u) +2\e(p-a_\infty)+2\e+\frac12\sum_{i=1}^2 (A_{ii}-B_{ii})\e_n\\
\neg&\neg\leq\neg &\neg\neg \l(\hat{w}) \int_{0}^{\hat{x}+\delta}[M_{\e}+2\e (\hat{x}-u+2\delta)]dG(u)+2\e( p-a_{\infty})+2\e+\frac12\sum_{i=1}^2 (A_{ii}-B_{ii})\e_n.
	\end{eqnarray*}
where $a_\infty:=\lim _{m\to +\infty}a_{\vsi_{m}}$. This is a contradiction when $\e$ and $\e_n$ are sufficiently small,
as $c>0$, $\int_0^{\hat x+\d} dG(u)\le 1$, and $M_\e\ge \vartheta>0$. That is, $\underline{u}\le \bar{u}^{\rho, \vth}$
must hold on $\sD_\d$.
%
%
The rest of the proof is straightforward, we leave it to the interested reader.
\qed

\begin{rem}
\label{u*}
{\rm
We recall that in \cite{BaiMa17} we proved the existence and uniqueness of the constrained viscosity solution. But the proof of the existence
essentially based on verifying that the value function is the desired viscosity solution. This fact sometimes  causes logical confusion, since a ``practical" version of the value function is actually the solution to HJB equation. Thus
is it often desirable, especially when an optimal strategy is based on the value function, to be able to ``construct" a constrained viscosity solution to the original problem, which we now describe.

First note that by uniqueness we need only show that we can construct a constrained viscosity subsolution $u^{*}$. Similar to the viscosity solution
of class ($\P$), we consider the class of constrained viscosity solution $v$ to (\ref{HJB0}) such that (i) $v(T,x,w)=0$; (ii) $x\mapsto v(t,x,w)$ is
increasing, for $\th=(t,x,w)\in D$; and (iii)  $v(t,x,w)$ is bounded on $D$, and $-Q_2T\le v(\th)\le (2+Q_1)T$, $\th\in D$, for some $Q_1, Q_2>0$. We shall call such viscosity solutions of class ($Q$).

Now let $d_{\sD}(\th):=\inf_{\eta\in {\sD}}|\eta-\th|$ be the distance between $\th$ and the set ${\sD}$. One can easily check that the functions
$\bar{\Upsilon}(\th)=2d_{\sD}(\th)+Q_{1}(T-s)$ and $\underline{\Upsilon}(\th)=d_{\sD}(\th)-Q_{2}(T-s)$, $\th\in D$, where
\bea
\label{Q12}
Q_1&=&\max\{2+M,2(p+\mu T)\};\\
Q_2&=&\Big[c+\sup_{0\leq w \leq T}  \Big|\frac{f(w)}{\overline{F}(w)}\Big|\Big]T+1,\nonumber
\eea
are respectively the viscosity supersolution on $\sD$ and subsolutions on $\sD^*$ to  (\ref{HJB0}) of class ($Q$)  with constants $(Q_1, Q_2)$.  Furthermore, $\underline{\Upsilon}\leq \bar{\Upsilon}$ on  $D$.
Now let ${\sM}$ be the set of all viscosity subsolution $u$ of (\ref{HJB0}) on $\sD^*$ of  class ($Q$) such that $\underline{\Upsilon}\leq u
 \leq \bar{\Upsilon}$, and define
$\frak{u}(s,x,w):= \sup_{u\in {\sM}} u(s,x,w)$. Then  similar to Theorem \ref{060501} one can show that $\frak{u}^{*}$ defined by
\bea
\label{fraku}
\frak{u}^{*}(s,x,w)=\limsup_{r\da 0}\{\frak{u}(t,y,v); (t,y,v)\in D \  \mbox{and}\  \sqrt{|t-s|+|y-x|^{2}+|v-w|^{2}} \leq r  \},
\eea
 is a (constrained) viscosity subsolution of (\ref{HJB0}) on $\sD^*$, and is of class ($Q$). In particular, by uniqueness (cf. \cite{BaiMa17}), $\frak{u}^*=V$,
 the value function of the original optimal dividend problem.
 \qed}
 \end{rem}

\section{The Regularity and Convergence of $\{V^{n,\delta}\}$. }
\setcounter{equation}{0}

We now turn our attention to the 
family $\{V^{n,\d}\}_{{n\ge 1,\d>0}}$, the solutions to the auxiliary equations (\ref{HJB2}). We shall argue that each $V^{n,\d}$ has desired the
regularity, and  $V^{n, \d}\to V$, the original value function in a satisfactory way, as $n\to\infty$ and $\d\to0$.

We first look at the regularity issue. To begin with, we note that if $u$ is a viscosity solution of (\ref{HJB2}) on $D_{\delta}$, and we consider
the change of variable: $y:=\ln(1+x+\d)$, $x\ge -\d$, and define $v(s,y,w):=u(s,e^y-1-\d,w)$, then it is easy to verify that $v$ is viscosity solution of
the  PDE:
\bea
\label{HJB2EE}
\dis v_{t}(\th)+\sup_{\gamma\in [0,1], a\in [0,M]}\sG^{n}(\th,v,v_y,v_w, v_{yy},v_{ww},I^{\delta}[v], \g, a)=0, \qq \mbox{ on $B_{\delta}$, }
\eea
where $\th=(s,y,w)$, $B_{\delta}:=\{\th=(s,y,w): 0\leq s<T+\delta, y>0, -\delta<w <s+\delta\}$, and
\bea
\label{Gn}
&&\sG^{n}(\th,v,v_y,v_w, v_{yy},v_{ww},I^\d[v], \g, a):=\Big[\frac{\e_{n}e^{-2y}}{2}+\frac{\sigma^{2}\gamma^{2}}2 \Big(\frac{e^{y}-\delta-1}{e^{y}}\Big)^{2}\Big]v_{yy}(\th)
+\frac{\e_{n}}{2}v_{ww}(\th)\nonumber\\
&&\qquad+\Big[pe^{-y}-\frac{\e_{n}}{2}e^{-2y}   -\frac{\sigma^{2}\gamma^{2}}2(\frac{e^{y}-\delta-1}{e^{y}})^{2}+(r+(\mu-r) \gamma) \frac{e^{y}-\delta-1}{e^{y}}\Big] v_{y}(\th)\\
 &&\qquad+a(1-e^{-y}v_{y}(\th))+v_{w}(\th)-cv(\th)+  \frac{f(w)}{\overline{F}(w)}I^{\delta}[v].  \nonumber
 \eea

It is worth noting that the main difference between (\ref{HJB2EE}) and (\ref{HJB2}) is that all the coefficients of (\ref{HJB2EE}) are bounded and continuous, and for each fixed $n\ge 1$ and $\d>0$, the function $\sG^n$ is uniformly {\it elliptic}. Therefore, a straightforward application of a combination of \cite[Lemma 2.9, Corollary 2.12 and Theorem 9.1]{CKS}
(see also \cite{W1} and  \cite[Theorem 1.1]{W2}) lead to the following result.
\begin{thm}
\label{regularity}
Assume Assumption \ref{assump1}. Let $u$ be the unique viscosity solution of class {\rm ($\tilde\P$)} to (\ref{HJB2EE}) with $\tilde\P(s,y,w):=\P(s,e^y-1-\d,w)$, $(s,y,w)\in D_\d$.
 Then,  $u\in \hC_{loc}^{2+\alpha}(D_{\delta})$\footnote{A function $u\in \hC^{1+\a}_{loc}([0,T]\times \hR)$ means $u\in L^\infty([0,T]\times \hR)$ and $Du\in \hC^\a_{loc}([0,T]\times \hR)$; $u\in \hC^{2+\a}_{loc}([0,T]\times \hR)$ means  $Du\in \hC^{1+\a}_{loc}([0,T]\times \hR)$, and $u_t$, $D^2u\in \hC^\a_{loc}([0,T]\times \hR)$.}
 in the sense that for any  compact set
$D^{'}\subset\neg\subset D_{\delta}$, there exists a constant $C>0$ such that
$\|u\|_{C^{2+\alpha}(D^{'})}\leq C$,
where $C>0$ depends on the uniform constants in Assumption \ref{assump1} and the time duration $T>0$.
\qed
	\end{thm}
\begin{rem}
\label{reguVnd}
{\rm
A direct consequence of Theorem \ref{regularity} is that the unique viscosity solution $V^{n,\d}$ to the PDE (\ref{HJB2}) in Theorem \ref{comp2}
has the same regularity for each fixed $n\ge 1$ and $\d>0$. This fact will be important for the construction of $\e$-optimal control in the sections
to follow.
\qed}
\end{rem}

 In the rest of the section we shall focus on an important and more involved issue: the convergence of the family $\{V^{n,\d}\}$, as $n\to \infty$ and
 $\d\to 0$. We shall first look at the limit as $n\to\infty$ (or as $\e_n\to0$). Naturally, let us consider an intermediate PDE:
   \bea
\label{HJB3}
\dis V_{t}(\th)+\sup_{\gamma\in [0,1], a\in [0,M]}H(\th,V, \nabla V, V_{xx}, V_{ww},{ I^{\d}[V]}, \g, a)=0,  \q  \th\in  \sD_{\delta}
\eea
where $H$ is defined by (\ref{H0}).
Following the same argument as that in \S2, we now argue that (\ref{HJB3}) admits a unique viscosity solution of class ($\P$). To see this,
for any $(t,y,v)\in D_{\delta}$, let
$$ \tilde{V}_{\delta}(t,y,v):=\lim_{k\rightarrow\infty} {\rm sup}\{V^{n,\delta}(\th):n\geq k, \th\in \bar{B}_{1/k}(t,y,v)\cap \bar{\sD}_{\d}\}, ~\mbox{\rm and}
$$
$$ \tilde{V}^{\delta}(t,y,v):=\lim_{k\rightarrow\infty} {\rm inf}\{V^{n,\delta}(\th):n\geq k, \th\in \bar{B}_{1/k}(t,y,v)\cap \bar{\sD}_{\d}\},$$
where  $B_r(t,y,v)$ is the open ball with radius $r$ centered at $(t,y,v)$, and $V^{n,\d}$'s are the viscosity solutions of class {\rm($\P$)} to PDE (\ref{HJB2}).
\begin{lem}
\label{120601}
For any $\Psi$ satisfying Assumption \ref{assump1}, the function $\tilde{V}_{\delta}$ (resp. $\tilde{V}^{\delta}$) is a viscosity subsolution (resp. supersolution) of class {\rm($\P$)} on $\sD_{\delta}$ to (\ref{HJB3}).
\end{lem}

{\it Proof.} We shall discuss only $\tilde{V}_{\delta}$ as the proof for $\tilde{V}^{\delta}$ is similar. First, it is easy to see that
$\tilde{V}_{\delta}$ is of class ($\P$) since all $V^{n,\d}$'s are uniformly bounded, uniformly in $n,\d$. Next, suppose that for some $\th_0:=(t_{0},y_{0},v_{0})\in D_{\delta}$,
 $0=[\tilde{V}_{\d}-\varphi](\th_0)$ is a (strict) maximum of $\tilde V_\d-\vf$
over $D_{\delta}$, where $\varphi\in C^{1,2,2}(\sD_{\delta})$.

For any $N>y_{0}$ we define $D_{\delta,N}=[0,T+\delta]\times [-\delta, N]\times [-\delta, s+\delta]$ so that $\th_0\in D_{\d, N}$. Since $\th_0$ is
the strict maximum of $\tilde V_\d-\vf$, for  $\e>0$, there exists a modulus of continuity $\o_{1}(\cd)$ such that
$$\sup_{\th\in B_{\varepsilon}^{c}(\th_0)\cap D_{\delta,N} }(\tilde{V}_{\delta}(\th)-\varphi(\th) )\leq -\o_{1}(\varepsilon)<0. $$
Now for $\bar\th:=(t,y,v)\in D_{\delta,N}$, by definition of $\tilde{V}_{\delta}$,  there exists $k_0:=k_{0}(\bar\th)=k_{0}(\bar\th;\e)$, such that
$$\sup_{ \th\in
\bar{B}_{1/k_{0}}(\bar\th)\cap \bar{\sD}_{\delta}}V^{n,\d}(\th)-\tilde{V}_{\delta}(\bar\th)<\frac{\o_{1}(\e)}{4} \qq n\geq k_{0}.
$$

Let us denote $\o_\vf^{\d, N}(\cd)$ to be
the modulus of continuity of $\vf$ on $D_{\delta,N}$. Then, for $\e>0$, there exists
$\eta_{0}:=\eta_{0}(\varepsilon)>0$ such that $\o^{\d, N}_\vf(\eta_0)<\o_1(\e)/4$. Thus, for $\bar\th\in D_{\delta,N}\backslash B_{\e}(\th_0)$ and  $n\geq k_{0}(\bar\th)$,
\beaa
&&\sup_{ \th
\in \bar{B}_{\frac1{k_{0}}\wedge\eta_{0}}(\bar\th)\cap \bar{\sD}_{\delta}}(V^{n,\delta}(\th)-\varphi(\th))=\sup_{{ \th
\in \bar{B}_{\frac1{k_{0}}\wedge\eta_{0}}(\bar\th)\cap \bar{\sD}_{\delta}}}(V^{n,\delta}(\th)-\tilde{V}_{\delta}(\bar\th)+\tilde{V}_{\d}(\bar\th)-\vf(\bar\th){+\vf(\bar\th)}-\vf(\th))\\
&&\qq\qq\qq \leq\frac{\o_{1}(\e)}{4}-\o_{1}(\e)+\o^{\d, N}_\vf({\eta_{0}})\leq\frac{\o_{1}(\e)}{4}-\o_{1}(\e)+\frac{\o_{1}(\e)}{4}=-\frac{\o_{1}(\e)}{2}.
\eeaa
 Since $B^c_{\e}(\th_0)\cap D_{\delta,N}$ is  compact and  $\bigcup_{\bar\th\in  {D}_{\delta,N}} B_{\frac1{k_{0}(\bar\th)}\wedge \eta_{0}}(\bar\th)\supset B^c_{\e}(\th_0)\cap D_{\delta,N}$,
there exists $N_{1}>0$  and $\th_{i}\in B^c_{\e}(\th_{0})\cap {D}_{\delta,N}$, $i=1,2,3...N_{1}$,
such that $\bigcup_{i=1}^{N_1} B_{\frac1{k_{0}(\th_i)}\wedge \eta_{0}}(\th_i)
 \supset  B^c_{\e}(\th_0)\cap D_{\delta,N}$.
 Hence, for any $n\geq \max_{1\leq i \leq N_{1}}k_0(\th_{i})$,
$$V^{n,\delta}(\bar\th)-\varphi(\bar\th)\leq -\frac{\o_{1}(\e)}{2}, \qq \bar\th\in B^c_{\e}(\th_0)\cap D_{\delta,N}.$$

Finally, let $\{\e_{\ell}\}_{\ell\in \hN}$ be a positive sequence  such that  $\e_{\ell}\downarrow 0$ as $\ell\rightarrow \infty$. For each
$\ell>0$, let $\bar\th_\ell \in B^c_{\e_1}(\th_0)\cap D_{\delta,N}$  and $n_{\ell}\geq \max\{{\max}_{1\leq i \leq N_{1}(\e_{\ell})}k_0(\th_{i}(\e_{\ell})),\frac{1}{\e_{\ell}}\} $ be
such that
\bea
\label{013101}
V^{n_{\ell},\d}(\bar\th_{\ell})-\vf(\bar\th_{\ell})=\max_{\bar\th\in \bar{D}_{\delta}}(V^{n_{\ell},\delta}(\bar\th)-\vf(\bar\th))>-\frac{\o_{1}({\e_{\ell}})}{2}.
\eea
Next, denoting $\vf^{n_{\ell}, \d}(\th):=\vf(\th)+V^{n_{\ell}, \d}(\bar\th_{\ell})-\vf(\bar\th_{\ell})$, $\th\in D_\d$, we see that $\vf^{n_{\ell}, \d}\in
\hC^{1,2,2}(D_{\delta})$, and
$0=V^{n_{\ell},\d}(\bar\th_{\ell})-\vf^{n_{\ell}, \d}(\bar\th_{\ell})=\max_{\th\in D_\d} V^{n_{\ell},\d}(\th)-\vf^{n_{\ell}, \d}(\th)$, and therefore
\bea
\label{013102}
&& \dis \varphi_{t}(\bar\th_{\ell})+\sup_{\g\in [0,1], a\in [0,M]}H^{n_{\ell}}(\bar\th_{\ell},\vf^{n_{\ell},\delta},
\nabla\vf, \vf_{xx},\vf_{ww},I^{\delta}[\vf^{n_{l},\delta}],\g, a) \ge 0.
\eea
Letting $\ell\rightarrow \infty$ in (\ref{013101}) and (\ref{013102}), we have
\beaa
0&\leq& \lim_{n_{\ell}\rightarrow \infty}V^{n_{\ell},\delta}(\bar\th_{\ell})-\varphi(\th_{0})\le \lim_{\e_{\ell}\to 0}\sup\{V^{n,\delta}(s,x,w):n\geq \frac1{\e_\ell}, (s, x,w)\in \bar{B}_{\e_{\ell}}(\th_{0})\cap \bar{\sD}_{\d}\}-\varphi(\th_{0})\\
&=& \lim_{k\rightarrow\infty}\sup\{V^{n,\delta}(s,x,w):n\geq k, (s, x,w)\in \bar{B}_{\frac1k}(\th_{0})\cap \bar{\sD}_{\d}\}-\varphi(\th_{0})=\tilde{V_{\delta}}(\th_{0})-\varphi(\th_{0})=0,
\eeaa
and
$\dis \varphi_{t}(\th_{0})+\sup_{\gamma\in [0,1], a\in [0,M]}H(\th_{0},\varphi,\nabla\varphi, \varphi_{xx},\varphi_{ww},I^\d[\varphi], \g, a)\geq 0$. That is, $\tilde V_\d$ is a viscosity subsolution of (\ref{HJB3}).
\qed

We should note that Lemma \ref{120601} and the comparison principle (Theorem \ref{comp2}) imply that $ \tilde{V}_{\delta}\leq \tilde{V}^{\delta} $.
On the other hand, by definitions of $ \tilde{V}_{\delta}$ and $\tilde{V}^{\delta} $, we also have $ \tilde{V}_{\delta}\geq\tilde{V}^{\delta} $. Thus
we have $ \tilde{V}_{\delta}=\tilde{V}^{\delta}$, and we shall denote it by $V^{\delta}$. Clearly, $V^\d\in \hC(D_{\delta})$.

Next, we recall the value function $V$ defined by (\ref{V1}). We know from \cite{BaiMa17} that it is the unique constrained
viscosity solution of (\ref{HJB0}), and from Remark \ref{u*} we see that it can be constructed as $\frak{u}^*$ defined by (\ref{fraku}).
In what follows we shall assume that, modulo a further approximation, we can always find a function $\Psi$ satisfying Assumption \ref{assump1},
such that $\Psi(\th)=\frak{u}^{*}(\th)=V(\th)$,  $\th\in \partial D$.
We should note that if $\Psi$ satisfies Assumption \ref{assump1}, then $\Psi$ will be smooth and having $\pa_x\Psi>1$ on the boundary $\pa D$.
However, these two conditions are {\it not} necessarily satisfied by the value function $V$. The following lemma is thus useful for our discussion.
\begin{lem}
\label{P=V}
Let $V$ be the value function defined by (\ref{V1}). Then there exists a
sequence of functions $\{\Psi_m\}_{m\ge 1}$ satisfying Assumption \ref{assump1}, and
continuous viscosity solutions $v^m$ of
 \bea
\left\{\ba{lll}
v_{t}(s,x,w)+\sL[v](s,x,w)=0, \qq
& (s, x,w)\in \sD,\\
v(s,x,w)=\Psi_{m}(s,x,w), & (s,x,w)\in \partial D.
\ea\right.
\eea
such that

(i) $\lim_{m\rightarrow\infty}\sup_{\th\in \partial D}|\Psi_{m}(\th)-V(\th)|=0$; and

(ii)  $\lim_{m\to\infty}\|v^{m}-V\|_{L^\infty(D)}\to 0$.
\end{lem}

{\it Proof. }  Let $V$  be the (viscosity) solution to (\ref{HJB0}) and $\vf_m:D\mapsto \hR$ the standard mollifiers of $V$. Then, since $V$ is continuous, we have
$\lim_{m\rightarrow\infty}\|\vf_{m}-V\|_{L^\infty(D)}=0$. Next, we define
  \bea
  \label{070502}
\Psi_{m}(\th)=\vf_{m}(\th)+(2+N_{m})d(\th,\partial D_{m}), \ \ \ \ \th:=(s,x,w)\in D
\eea
where $N_{m}:=\sup_{(s,w)\in[0,T]\times[0,s]}|\pa_x\vf_{m}(s,0,w)|$, and $\{D_{m}\}_{m\ge 1}$ is a sequence of smooth area such that $D\subset D_{m}$,
$d(D,D_{m})<\delta_{m}:=\frac{1}{m(2+N_{m})}$, and $D_{m}$ is parallel to the plane $\{(s,x,w),-\delta_{m}\leq s \leq T+\delta_{m}, x=0, -\delta_{m}\leq w \leq s+\delta_{m}\}$. It is then easy to check that $\sup_{\th\in \partial D}|(2+N_{m})d(\th,\partial D_{m})|\le \frac2m$, and
$\pa_x\Psi_{m}(s,0,w)= \pa_xf_{m}(s,0,w)+(2+N_{m})\ge -N_{m}+2+N_{m}=2$. Consequently, one can further check that,
by defining $\Psi_m\equiv0$ on $D^c_1$, all $\Psi_m$'s satisfy Assumption \ref{assump1}.
Now let   $v^{m}$  be the unique viscosity solution of (\ref{HJB0}) on $\sD$ with $v^m=\Psi_m$ on $\pa D$.
Then by definition (\ref{070502}) we can easily check that $a_{m}:=\sup_{\th\in \partial D}|v^{m}(\th)-V(\th)|=\sup_{\th\in \partial D}|\Psi_{m}(\th)-V(\th)|\to 0$, as $m\to\infty$, and
$ v^{m}-a_{m}\leq V \leq v^{m}+a_{m}$, on $\partial D$.
Since $v^{m}-a_{m}$ and $v^{m}+a_{m}$ are the viscosity subsolution and  supersolution of (\ref{HJB0}) on $\sD$, respectively,  by
comparison theorem we can then deduce that $\lim_{m\to\infty}\|v^{m}-V\|_{L^\infty(D)}\to 0$, proving the lemma.
\qed

We can now prove the main result of this section.
\begin{thm}
\label{convg}
Let $V$ be the value function defined by (\ref{V1}). Then
for any $\e>0$, there exists $n\in\hN$, and $\d>0$, depending only on $\e$, such
that
$\|V^{n,\d}-V\|_{L^\infty(D)}<\e$,
where $V^{n, \d}\in \hC^{2+\a}(D_\d)$ is a (viscosity) solution to (\ref{HJB2}) of class {\rm ($\P$)}, for some function $\Psi$ satisfying Assumption \ref{assump1}.
\end{thm}

{\it Proof.} In light of Lemma \ref{P=V},  we can assume without loss of generality that we can find $\Psi$ satisfying Assumption
\ref{assump1} such that $\Psi=\frak{u}^*=V$ on $\pa D$. (Otherwise for any $\e>0$ we can first choose $\Psi_m$ so that it satisfies Assumption
\ref{assump1}, and the corresponding viscosity solution $v^m$ satisfies $\Psi_m=v^m$ on $\pa D$, and $\|v^m-V\|_{L^\infty(D)}<\e/3$, and then prove the theorem for $\Psi_m$ and $v^m$.)
For convenience we shall also define
$\frak{u}^{*}(\th)=\Psi(\th)$ for $\th\in \sD^{*,c}_\d$ (see (\ref{Ddc})).

Now let $V^{n,\delta}$ be the solutions of (\ref{HJB2}) of class ($\P$).
We first show that  $\lim_{n\to\infty}\|V^{n,\delta}-V^{\delta}\|_{L^{\infty}(D_{\delta})}= 0$.
 Indeed, if not, then there exist  $\varepsilon_{0}>0$,
$\{n_{k}\}_{k\in \hN}\subset \hN$, and  $\{\th_k:=(t_{k},x_{k},w_{k})\}_{k\in \hN}\subset D_{\delta}$, such that $n_{k}\uparrow\infty$, as $k\rightarrow \infty$, and
$$|V^{n_{k},\delta}(t_{k},x_{k},w_{k})-V^{\delta}(t_{k},x_{k},w_{k})|>\varepsilon_{0}.$$
By definition of $\bar D_\d$ we see that, taking a subsequence if necessary, we can assume that there exists $\th_0:=(t_{0},x_{0},w_{0})\in \bar{D_{\delta}}$ (allowing $x_0=+\infty$) such that $\th_{k}\rightarrow \th_{0}$. Now let $k\rightarrow \infty$.
If $x_{0}<+\infty$, the we have
$\tilde{V}_{\delta}(\th_{0})-V^{\delta}(\th_{0})\geq \varepsilon_{0}$ or $\tilde{V}^{\delta}(\th_{0})-V^{\delta}(\th_{0})\leq-\varepsilon_{0}$, which contradicts the fact that
$\tilde{V}_{\delta}=\tilde{V}^{\delta}=V^{\delta}$ in $D_{\delta}$. If $x_{0}=+\infty$, then we have
$\tilde{V}_{\delta}(t_{0},N,w_{0})-V^{\delta}(t_{0},N,w_{0})\geq \varepsilon_{0}$ or $\tilde{V}^{\delta}(t_{0},N,w_{0})-V^{\delta}(t_{0},N,w_{0})\leq-\varepsilon_{0}$, for some $N>0$, also a contradiction. This proves the claim.

 Next, let us denote $a_{\delta}:=\sup_{\th\in D_{\delta }\backslash D}|V^{\delta}(\th)-V(\th)|$. Then, noting that $\P=V=
\frak{u}$ on
$\pa D_\d$,
for $\bar\th=(t,y,v)\in \pa D_\d$, we have
\beaa
a_{\d}
&=& \sup_{\th\in D_\d\setminus D}|V^\d(\th)-\psi(\bar\th)+\psi(\bar\th)-V(\th)|
\le \sup_{\th\in D_\d\setminus D}[|V^{\d}(\th)-V^\d(\bar\th)|+|\psi(\bar\th)-\psi(\th)|]\\
&\leq& \sup_{\th\in D_\d\setminus D}[\o(|\th-\bar\th|)+|\psi(\bar\th)-\psi(\th)|]=o_{\delta}(1), \q\mbox{ as $\delta\rightarrow 0$.}
\eeaa
Here $\o(\cd)$ is the modulus of continuity of $V^{n, \d}$ (which can be chosen to be independent of $\d$(!)).
Furthermore, it is easy to verify that $V^{\delta}-a_{\delta}$ and $V^{\delta}+a_{\delta}$ are  viscosity subsolution and viscosity supersolution of (\ref{HJB3}), respectively,  and $V^{\delta}-a_{\delta}\leq V \leq V^{\delta}+a_{\delta}$, on $\partial D$. It then follows from the comparison principle that $\|V^{\delta}-V\|_{L^{\infty}(D)}=o_{\delta}(1)$, as $ \delta\rightarrow 0$.

Combining above, for $\e>0$, we can first choose $\d=\d(\e)>0$ so that $\|V^{\delta}-V\|_{L^{\infty}(D)}<\e/2$, and then choose $n=n(\d(\e))\in \hN$ such that $\|V^{n,\delta}-V^{\delta}\|_{L^{\infty}(D)}\le \|V^{n,\delta}-V^{\delta}\|_{L^{\infty}(D_{\delta})}<\e/2$. Noting that $V^{n, \d}\in \hC^{2+\a}_{loc}(D_\d)$, thanks to  Theorem \ref{regularity} and Remark \ref{reguVnd}. The proof is now complete.
\qed

\section{{ Construction of $\varepsilon$-Optimal Strategy}}
\setcounter{equation}{0}

We are now ready to construct the desired $\e$-optimal strategy. The idea is simple: for each $\e>0$, we choose an
approximating solution $V^{n, \d}$, guaranteed by Theorem \ref{convg}, and define a strategy in the form of (\ref{optim}). It is then reasonable
to believe that such a strategy should be $\e$-optimal.

To be more precise, let  $\{\e_k\}$ be any sequence such that $\e_k\da 0$, as $k\to \infty$, and let $V^k:=V^{n_{k},\d_{k}}\in \hC^{2}_{loc}(D_{\d_k})$ be the corresponding solutions of (\ref{HJB2}) as those in Thoerem \ref{convg}. That is,
%
\bea
\label{Vk}
\|{V}^{n_k,\d_k}-V\|_{L^{\infty}(D)}<\e_k \rightarrow 0, \qq \mbox{as $k\to\infty$.}
\eea
Since $V(\th)\equiv 0$ for $\th\in D^c$, we can and shall assume that $V^k(\th)\equiv0$ for $\th\in D^c$, for all $k$.  Furthermore, since each
 ${V}^{n,\delta}$ is of class ($\Psi$) for some $\Psi$ satisfying Assumption \ref{assump1}, we can assume  ${V}^{n,\delta}_{x+}(s,-\delta,w)> 1$.
Therefore ${V}^{k}_{x+}(s,0,w)> 1$ for large $k$.
%

Now recall the optimal strategy (\ref{optim}). We consider the sequence of strategies $\{(\g^k, a^k)\}_{k\in\hN}$:
\bea
\label{gkak}
\left\{\ba{lll}
\dis \g^k_t:=\frac{(\mu-r)V^{k}_{x}(t, X_t, W_t)}{\sigma^{2}X_t V^{k}_{xx}(t, X_t, W_t)}; \ms\\
a^k_t:= M \b1_{\{V^{k}_{x}(t, X_t, W_t)<1\}}+p \b1_{\{V^{k}_{x}(t, X_t, W_t)=1\}},
\ea\right.
\eea
where $(X,W)$ is the solution to the corresponding close-loop system (\ref{Xsxw}). More precisely, let us define, for each $k\in\hN$,
two functions $\G^k,\Xi^k: [0,T]\times\hR\times \hR\mapsto\hR$:
\bea
\label{GaXi}
\G^k(s,x,w):=\frac{(\mu-r)V^{k}_{x}(s,x,w)}{\sigma^{2}x V^{k}_{xx}(s,x,w)}; \q \Xi^k(t,x,w):= M \b1_{\{V^{k}_{x}(s, x,w)<1\}}+p \b1_{\{V^{k}_{x}(s,x, w)=1\}}.
\eea
Then $(\g^k_t, a^k_t)=(\G^k(\Th^k_t), \Xi^k(\Th^k_t))$, $t\in [0,T]$, where $\Th^k_t=(t, X^k_t,W_t)$, and $(X^k, W)$ is the, say, weak solution to the
``close-loop" dynamics of the reserve (recall (\ref{Xsxw})), defined on some probability space $(\O, \cF, \hP, \hF)$:
\bea
\label{closeloop}
\left\{\ba{lll}
dX_t=b^k(t, X_t, W_t)dt+\si^k(t, X_t, W_t)dB_t-dQ^{s,w}_t, \q X_s=x; \ms\\
W_t=w+(t-s)-(\si_{N_t}-\si_{N_s}),  \qq 0\le s\le t\le T.
\ea\right.
\eea
Here in the above,  for $\th:=(s,x,w)\in D$, we have (noting that ${V}^{k}_{x}(s,x,w)>0$ for $(s,x,w)\in D$)
\bea
\label{approxstratb}
b^k(\th)&:=&\left\{\ba{lll}
p+rx-(\mu-r)\G^k(\th)x-\Xi^k(\th) \qq\qq &0<\G^k(\th)\leq 1\ms\\
p+\mu x-\Xi^k(\th),  &\mbox{otherwise;}
\ea\right.
\ms
\\
\label{approxstrata}
\si^k(\th)&:=&\left\{\ba{lll}
\si x\G^k(\th) \qq &0<\G^k(\th)\leq 1\ms\\
\sigma x,  &\mbox{otherwise.}
\ea\right.
\eea

In the rest of the section we shall verify the main result of this paper: (i)  the close-loop system (\ref{closeloop}) is well-posed, and (ii)
$(\g^k, a^k)$ provides an $\e$-optimal  strategy for $k$ large.
 To this end, we make some quick observations on the functions $\G^k$ and $\Xi^k$ in (\ref{GaXi}) that determines the optimal strategy. Clearly,
 the function $\G^k$ is  continuous,
 and $\G^k(\th)=1$, when $x$ is close to 0. In fact, by a further approximation (cf. \cite{LSM}) if necessary, we can actually assume further that $\G^k$ is Lipschitz continuous (with Lipschitz constant depending on $k$).  The main difficulty, however, is that the function $\Xi^k$ is {\it discontinuous}. This, together with the presence of jumps, makes finding the strong solution to SDE (\ref{closeloop}) a much more involved task.
Our plan of attack is the following. We shall begin by looking at the {\it weak solution} to (\ref{closeloop}). Then using the fact that the SDE is
one-dimensional, we shall argue that the weak solution is actually strong as well as pathwisely unique, up to the ruin time $\t= \inf\{t>0, X_t<0\}$.

To do this, we shall modify the function $\si^k$ slightly: for $m\in\hN$, we consider $\vf^m(x)=\frac1{m}\vee x\wedge m$, and define
$\si^{m, k}(\th):=\si\vf^m(x)\G^k(\th)$, $\th\in D$. Since both $\vf^m$ and $\G^k$ are bounded and Lipschitz, so is
$\si^{m,k}$. Furthermore, it is readily seen that for some constant $c_m>0$, one has
 \bea
 \label{sik}
0< c_m\le \si^{m,k}(\th)\le \si (x\wedge m)  \q\th:=(s,x,w)\in D.
 \eea

To continue our discussions let us now consider the {\it canonical space}. Let $\O^1=\hC([0,T])$, the space of all continuous functions, null at zero, and endowed with the usual sup-norm.
Let $\cF_t^{1}\dfnn\si\{\o(\cd\wedge t)|~\o\in\O^{1}\}$, $t\ge 0$, $\cF^{1}\dfnn\cF_T^{1}$, $\hF^1=\{\cF^{1}_{t}\}_{t\in [0,T]}$ and $\hP^0$ be the Wiener measure on $(\O^{1}, \cF^{1})$ so that the {\it canonical process}
$ B_t(\o)\dfnn\o^1(t)$, $(t,\o^1)\in[0,T]\times \O^{1}$ is an $(\hP^0, \hF^1)$-Brownian motion.
Let $\O^{2}=\hD([0,T])$, the space of all real-valued, c\`adl\`ag (right-continuous with left limit) functions, endowed
with the Skorohod topology, and  similarly define
$\hF^2=\{\cF^{2}_{t}\}_{t\in [0,T]}$ and $\cF^{2}\dfnn\cF_T^{2}$.
Let $\hP^{Q}$ be the law of the renewal claim process $Q$ on $\hD([0,T])$, so that the  coordinate process,
$Q_t(\o^2)=\o^2(t)$, $(t,\o^2)\in[0,T]\times\O^{2}$.
Now we consider the product space:
\bea
\label{WP}
\O\dfnn\O^1\times \O^2;\q \cF\dfnn\cF^1\otimes\cF^2;\q \hP\dfnn\hP^{0}\otimes\hP^{Q}; \q \cF_t\dfnn\cF^1_t\otimes\cF^2_t, ~~t\in[0,T].
\eea
We now consider the following SDE on the canonical space $(\O^1,\cF^1,\hP^0;\hF^1)$:
\bea
\label{SDE2}
\left\{\ba{lll}
\dis dX_{t}=\si^{m,k}(t,X_{t}, W_t)dB_t-dQ_t, \qq X_0=x;\ms\\
W_t=t-\si_{N_t},
\ea\right. \qq   t\in[0,T].
 \eea
We have  the following result.
\begin{prop}
\label{Xmk}
Assume Assumption \ref{assump0}. Then, the SDE (\ref{SDE2}) has a strong solution.
\end{prop}

{\it Proof.} We write the element of $\O$ as $\o=(\o^1,\o^2)\in\O$. Then, the two marginal coordinate processes defined by
$B_t(\o)\dfnn \o^1(t)$, $ Q_t(\o)\dfnn \o^2(t)$, $(t,\o)\times [0,T]\times\O$. Then under our assumptions $B$ and $Q$ are independent, and the process $Q_t(\o)=\o^2(t)$ is piecewise constant jumping
 at $0<\sigma_{1}(\o^2)< \cdots<\sigma_{N_T(\o^2)}(\o^2)<T$, where $N_t(\o^2)$ denotes the number of jumps of $Q$ up
to time $t$, hence a renewal counting process. We then  define  $W_t(\o)=t-\si_{N_t(\o^2)}(\o^2)$, $t\ge 0$.

Now on the canonical process, for $\hP^Q$-a.s. $\o^2\in \O^2$ we define
\bea
\label{bOx}
\tilde\si^{m,k,\o^2}(t,x):=\si^{m,k}(t, x-\o^2(t), t-\si_{N_t(\o^2)}(\o^2)), \qq (t,x)\in [0,T]\times \hR,
 \eea
and consider  the SDE on the space $(\O^1,\cF^1,\hP^0;\hF^1)$:
\bea
\label{SDE-0}
\tilde X_{t}=x+\int_{0}^{t}\tilde \si^{\o^2, m,k}(s,\tilde X_{s})dB_s, \qq   t\in[0,T]{\color{purple2}.}
 \eea
Clearly, by definition (\ref{bOx}) and the facts (\ref{sik}) and that $\si^{m,k}$ is Lipschitz, SDE (\ref{SDE-0}) has a unique strong solution $\tilde X^{\o^2}_t:=\tilde X_t(\cd, \o^2)$ on $(\O^1,\cF^1,\hP^0;\hF^1)$, for $\hP^Q$-a.s. $\o^2\in\O^2$. Consequently, by (\ref{bOx}), if we define $X:=\tilde X-Q$, and $W_t=t-\si_{N_t}$, then $( X, W)$ satisfies (\ref{SDE2}).

The uniqueness of the solution $(X,W)$ follows from that of $\tilde X$ as $Q$ is a coordinate process, completing the proof.
\qed

Now let  $(X, W)$ be a strong solution of (\ref{SDE2}) on $(\O, \cF, \hP)$, and denote it by $(X^{m,k}, W)$ if the dependence on $m,k$ is important. Clearly, for fixed $\o^2\in \O^2$, $X^{m,k}_t(\o)=\tilde X^{\o^2}_t-\o^2(t)$. It is well-known  (cf. e.g., \cite{Aronson1} and \cite{IKO}) that the solution $\tilde X^{\o^2}$ of (\ref{SDE-0}) has a
transition density, denoted by $p^{\o^2}(t,y;s,x)$ to indicate its dependence on $\o^2$, and it satisfies
\bea
\label{densityX}
p^{\o^2}(t,y;s,x)\leq M_0 |t-s|^{-\frac{1}{2}}{\rm exp}\{\frac{-\Lambda (y-x)^{2}}{t-s}\}, \q s\le t, ~x,y\in\hR,
\eea
where constants $M_0$ and $\Lambda$ depends only on $m, k$, but independent of $\o^2$. Consequently, for fixed $\o^2\in\O^2$, $X^{m,k}(\cd, \o^2)$ has the density function $p^{\o^2}(t, y+\o^2(t);s, x)$ under $\hP^0$. Furthermore, by
renewal theory (see, e.g., \cite{Ross}), the random variable $\sigma_{N_{t}}$ has a density function
\bea
\label{fsiN}
f_{\si_{N_{t}}}(u)=\bar F(t-u)m'(u)=\bar{F}(t-u)\sum_{n=1}^{\infty}f_{n}(u),\qq t\geq u \geq 0,
\eea
where $m(t)=\hE[N_t]=\sum_{n=1}^{\infty}F_{n}(t)$, $F$ is the law of the waiting time $T_i$'s, $F_n$ is the n-fold convolution of $F$ with itself,  and $f_{n}$ is corresponding density function.
Therefore, we can write down the joint distribution of
$(X^{m,k},\si_{N_t})$:
\bea
\label{disXW}
\hP(X^{m,k}_t\in A, \si_{N_t}\in B)&=&\int_{\O^1}\int_{\O^2}\1_{\{X^{m,k}_t(\o^1, \o^2)\in A\}}\1_{\{\si_{N_t}(\o^2)\in B\}}\hP^0(d\o^1)\hP^Q(d\o^2)\nonumber\\
&=&\int_{\O^2}\Big[\int_{A}p^{\o^2}(t, y+\o^2(t);s,x)dy\Big]\1_{\{\si_{N_t}(\o^2)\in B\}}\hP^Q(d\o^2).
\eea

In what follows we shall make use of an extra assumption on the jump times $\si_{N_t}$.
\begin{assum}
\label{assump2}
There exist constant $\g^{'}>1$ such that
 \bea
\label{jump}
\int_0^T\int_0^t t^{\frac{1-\g'}{2}}f^{\g'}_{\si_{N_t}}(u)dudt<+\infty.
\eea
\end{assum}
\begin{rem}
\label{remUi}
{\rm We remark that the Assumption \ref{assump2} is merely technical, but it covers a large class of cases that are commonly seen in applications. In particular, we note that if we take $\frac{3-\gamma^{'}}{2}>-1$, then $\gamma^{'}<5$. Furthermore,
if $T_{i}$ is of exponential distribution with $\lambda$ (that is, the renewal process $N$ becomes Poisson), then $m(t)=\hE N(t)=\lambda t$ and $
f_{\si_{N_{t}}}(u)=\lambda e^{-\lambda (t-u)}$. Then,
\beaa
\int_0^T\int_0^t t^{\frac{1-\g'}{2}}f^{\g'}_{\si_{N_t}}(u)dudt=\int_0^T\int_0^t t^{\frac{1-\g'}{2}}(\lambda e^{-\lambda (t-u)})^{\g'}dudt
\leq  \int_0^T\int_0^t t^{\frac{1-\g'}{2}}\lambda ^{\g'}dudt=\frac{2\lambda^{\g'}}{5-\gamma^{'}}T^{\frac{5-\g'}{2}}.
\eeaa

Also, if $T_{i}\sim$  Erlang(k,$\lambda$), that is, $F(u,k,\lambda)=1-\Sigma_{i=0}^{k-1}\frac{1}{i!}e^{-\lambda x}(\lambda x)^{i}$, as we
often see in the Sparre Andersen models, then $\sum_{n=1}^{\infty}f_{n}(u,k,\lambda)\leq \sum_{n=1}^{\infty}f_{n}(u,1,\lambda)=\lambda $, and
one can check that
\beaa
\int_0^T\int_0^t t^{\frac{1-\g'}{2}}f^{\g'}_{\si_{N_t}}(u)dudt\leq \frac{2\lambda^{\g'}}{5-\gamma^{'}}T^{\frac{5-\g'}{2}}.
\eeaa
In both cases Assumption \ref{assump2} holds.
\qed
}
\end{rem}

\section{Strong Well-posedness of the Closeloop System}
\setcounter{equation}{0}

We  now ready to study the existence and (pathwise) uniqueness of the close loop system (\ref{closeloop}). Again, for each $m\in\hN$ we consider the ``truncated" version
of $b^k$: $b^{m, k}(t,x,w):=\beta^k(t, -m\vee x\wedge m, w)$. Then $b^{m, k}$ is a bounded and measurable function.
Let $(X^{m, k}, W)$ be the strong solution of (\ref{SDE2}) on $(\O^1, \hF^1, \hP^0)$, and for  $\o^2\in\O^2$, define
\bea
\label{thmk}
\th^{m,k}_t(\cd, \o^2):=\frac{b^{m, k}(t, X^{m,k}_t(\cd, \o^2), t-\si_{N_t(\o^2)}(\o^2))}{\si^{m, k}(t, X^{m,k}_t(\cd, \o^2), t-\si_{N_t(\o^2)}(\o^2))}.
\eea
Since $b^{m,k}$ is bounded, by (\ref{sik}) we see that, modulo a $\hP^Q$-null set $N^2\subset \O^2$, $\th^{m,k}(\cd, \o^2)$ is an bounded, $\hF^1$-adapted process, for all
$\o^2\in \O^2\setminus N^2$. We can then define the following exponential martingale on $(\O^1, \cF^1, \hP^0;\hF^1)$:
\bea
\label{Mmk}
L^{m,k}_{t}(\cd, \o^2):=\exp\Big\{\int_{0}^{t}\th^{m,k}_s(\cd, \o^2)dB_{s}-\frac{1}{2}\int_{0}^{t}|\th^{m,k}_s(\cd, \o^2)|^2ds\Big\},
\qq \o^{2}\notin N^2,
\eea
and a new probability measure $\tilde{\hP}^{m,k}$ on $(\O, \cF)$ by
\bea
\label{tildeP}
\tilde \hP^{m,k}(A_1\times A_2):=\int_{A_2}\int_{A_1}L_{T}^{m,k}(\o^1,\o^2)\hP^0(d\o^1)\hP^Q(d\o^2), \q A_1\in\cF^1,~A_2\in \cF^2.
\eea
 Then, it is readily seen that,  under $\tilde{\hP}^{m,k}$, $\tilde{B}^{m,k}_t:=B_t-\int_{0}^t\th^{m,k}_sds$, $t\in[0,T]$, is a Brownian motion,  still independent of $Q$, and on the space $(\O, \cF, \tilde\hP^{m,k})$, $(X^{m,k}, W)$ satisfies, for $t\in [0,T]$,
 \bea
 \label{SDE4}
\left\{\ba{lll}
\dis dX^{m,k}_t=b^{m,k}(t, X^{m,k}_t, W_t)dt+\si^{m,k}(t, X^{m,k}_t,W_t)d\tilde B^{m,k}_t-dQ_t, ~ X^{m,k}_0=x\ms\\
W_t=t-\si_{N_t}.
\ea\right.
\eea
In other words, $(\O, \cF, \tilde\hP^{m,k}, \tilde B^{m,k}, X^{m,k}, W)$ is a {\it weak solution} to a truncated version of (\ref{closeloop}). Our task in this subsection is to show that
this weak solution can actually be strong and that it is pathwisely unique. Furthermore, we shall argue that, as $m\to\infty$, the sequence $\{X^{m,k}\}$ would
converge to a process $X^k$, which satisfies the SDE (\ref{closeloop}) on the interval $[0,\t_k)$, where $\t_k:=\inf\{t>0: X^k_t<0\}$. This is clearly sufficient for our purpose.

We should note that since the coefficient $b^{m,k}$ is
discontinuous, the pathwisely unique strong solution is only possible because the SDE (\ref{closeloop}) is one-dimensional. Our argument borrows the idea initiated in
Gy\"ongy and Pardoux \cite{GP} (see also, e.g., \cite{BM}), using the so-called Krylov estimate (cf. \cite{Krylov}). To this end,
let us begin with some observations. Let $(X^{m,k},W, B)$ be any weak solution of SDE (\ref{SDE4}) defined on some filtered probability space $(\O, \cF, \hP;\hF)$, we may assume that $(\O, \cF)$ is
the canonical space defined before, except that $\hP$ is any probability measure, and $\hF$ is augmented by all the $\hP$-null sets. Recalling  $\th$ and $M$ defined by (\ref{thmk}) and (\ref{Mmk}), respectively, define $\bar \th:=-\th$, and $\bar L:=L^{-1}$. Note
that the process $\th$ actually depends on $\o^2$, namely we should have $\th=\th^{\o^2}$, for $\o^2\in \O^2$, and hence
$L=L^{\o^2}$ as well. We now define, for fixed $\o^2$, a new probability measure $\frac{d\hP^{0,\o^2}}{d\hP}\big|_{\cF^1_T}=\bar L^{\o^2}_T$ on
$(\O^1, \cF^1)$, so that  $B^0_t:=B_t-\int_0^t\bar\th^{\o^2}_s ds$, $t\ge 0$ is a Brownian motion
on $(\O^1, \cF^1, \hP^{0,\o^2})$. We next define a new probability measure on $(\O, \cF)$ by
\bea
\label{barP}
\bar\hP(A\times B)=\int_B\int_A\hP^{0,\o^2}(d\o^1)\hP^Q(d\o^2)=\int_B\int_A\bar L^{\o^2}_T(\o^1)\hP(d\o^1\otimes d\o^2), ~ A\in \cF^1, B\in\cF^2.
\eea
Then it is readily seen that $\bar L_t(\o)=\bar L_t(\o^1, \o^2):=[L^{\o^2}]^{-1}(\o^1)$, $t\in[0,T]$ is a martingale under $\bar\hP$, $\frac{d\bar\hP}{d\hP}\big|_{\cF_T}=\bar L_T$, and $(X^{m,k}, W, B^0)$ solves SDE (\ref{SDE2}) on the  space $(\O, \cF, \bar\hP)$.

We are now ready to prove the following {\it Krylov estimate}.
\begin{lem}
\label{krylov1}
Assume Assumptions \ref{assump0} and \ref{assump2}. Let
$X^{m,k}$ be a weak solution of SDE (\ref{SDE4}). Then,
 for any bounded and measurable function $g:[0,T]\times[0,+\infty)\times [0,T]\rightarrow \hR_+$, it holds that
\bea
\label{krylovI}
\hE\int_{0}^{T}g(t,X_{t},W_{t})dt\leq G  \Big\{\int_{0}^{T}\int_{\hR}\int_{0}^{t}g^{\b\gamma}(t,y,t-u)dydu dt\Big\}^{1/\b\gamma}.
\eea
Here in the above $G$ is a constant defined by
\bea
\label{G}
G=C(M_0,\L,\g',\b)\{\bar\hE\bar L_{T}^{-\a}\}^{1/\a}\Big[\int_0^T\int_0^tt^{\frac{1-\g'}{2}}f^{\g'}_{\si_{N_t}}(u)dudt\Big]^{\frac1{\b\g'}},
\eea
where $\bar \hE=\hE^{\bar\hP}$, and $\bar L_T=\frac{d\bar\hP}{d\hP}$, $\frac1\a+\frac1\beta=1$, and $\g'$ is given in Assumption \ref{assump2}.
\end{lem}
{\it Proof.}  Throughout this proof we fix $m$ and $k$, and thus omit them in the notation for simplicity.
%
For any  bounded, nonnegative measurable function $g:[0,T]\times[0,+\infty)\times [0,T]\rightarrow \hR_+$ we have
 \bea
 \label{Kryest}
&&\hE\Big\{\int_{0}^{T}g(t,X_{t},W_{t})dt\Big\}=\bar\hE\Big\{\bar L^{-1}_T\int_{0}^{T}g(t,X_{t},W_{t})dt\Big\}\nonumber\\
&\leq &\{\bar\hE\bar L_{T}^{-\a}\}^{1/\a}\Big\{\bar\hE\Big[\int_0^T g^{\b}(t,X_t, t-\si_{N_t})dt\Big]\Big\}^{1/\b}\\
&= &\{\bar\hE\bar L_{T}^{-\a}\}^{1/\a}\Big\{\int_0^T\int_{\O^2}\Big[\int_{\hR}g^{\b}(t,y,t-\si_{N_t}(\o^2))p^{\o^2}(t,y+\o^2(t),0,x)dy\Big]\hP^Q(d\o^2)dt\Big\}^{1/\b}\nonumber\\
&\leq &\{\bar\hE\bar L_{T}^{-\a}\}^{1/a}\Big\{\int_0^T\int_{\O^2}\Big[\int_{\hR}g^{\b}(t,y,t-\si_{N_t}(\o^2))M_0 |t|^{-\frac{1}{2}}e^{\frac{-\Lambda (y+\o^2(t)-x)^{2}}{t}}dy\Big]\hP^Q(d\o^2)dt\Big\}^{1/\b}.\nonumber
\eea
Note that, by H\"older's inequality again, we have
\bea
\label{Kryest1}
&&\int_{\hR}g^{\b}(t,y,t-\si_{N_t}(\o^2))M_0 |t|^{-\frac{1}{2}}\exp\Big\{\frac{-\Lambda (y+\o^2(t)-x)^{2}}{t}\Big\}dy\\
&\leq &\Big[\int_{\hR}g^{\b\g}(t,y,t-\si_{N_t}(\o^2))dy\Big]^{\frac{1}{\g}}\Big[\int_{\hR}\Big(M_0 |t|^{-\frac{1}{2}}\exp\Big\{\frac{-\L(y+\o^2(t)-x)^{2}}{t}\Big\}\Big)^{\g'}dy\Big]^{\frac{1}{\g'}}\nonumber
\eea
where  $1/\gamma+1/\gamma^{'}=1$. By the direct calculation, we have
 \bea
 \label{Kryest2}
\int_{\hR}\Big(M_0 |t|^{-\frac{1}{2}}\exp\Big\{\frac{-\L(y+\o^2(t)-x)^{2}}{t}\Big\}\Big)^{\g'}dy\leq C(M_0,\Lambda,\gamma^{'})|t|^{\frac{1-\gamma^{'}}{2}},
 \eea
where $C(M_0, \L, \g')$ is some constant depending only on $M_0$, $\L$, and $\g'$. Plugging (\ref{Kryest1}) and (\ref{Kryest2}) into (\ref{Kryest}), and applying the H\"older inequality again we obtain that
\bea
\label{Kryest3}
&&\hE\Big\{\int_{0}^{T}g(t,X_{t},W_{t})dt\Big\} \nonumber\\
&\le& \{\bar\hE\bar L_{T}^{-\a}\}^{1/\a}\Big\{\int_0^T\int_{\O^2}\Big[\int_{\hR}g^{\b\g}(t,y,t-\si_{N_t}(\o^2))dy\Big]^{\frac{1}{\g}}C(M_0,\L,\g')|t|^{\frac{1-\g^{'}}{2\g'}}\hP^Q(d\o^2)dt\Big\}^{1/\b} \nonumber\\
&\leq & C(M_0,\L,\g',\b)\{\bar\hE\bar L_{T}^{-\a}\}^{1/\a}\Big\{\int_0^T\int_0^t\Big[\int_{\hR}g^{\b\gamma}(t,y,t-u)dy\Big]^{\frac{1}{\g}}
|t|^{\frac{1-\g'}{2\g'}}f_{\si_{N_t}}(u)dudt\Big\}^{\frac1\b}\\
&\le&G\Big[\int_0^T\int_0^t\int_{\hR}g^{\b\gamma}(t,y,t-u)dydudt\Big]^{\frac{1}{\b\g}}, \nonumber
\eea
where  $C(M_0,\L,\g',\b):= C^{\frac1\b}(M_0,\L,\g')$, and
\bea
\label{G}
G:=C(M_0,\L,\g',\b)\{\bar\hE\bar L_{T}^{-\a}\}^{1/\a}\Big[\int_0^T\int_0^tt^{\frac{1-\g'}{2}}f^{\g'}_{\si_{N_t}}(u)dudt\Big]^{\frac1{\b\g'}}.
\eea
The proof is now complete.
\qed
%

We are now ready to prove that, for fixed $m,k$, SDE (\ref{SDE4}) actually has a pathwisely unique strong solution on the interval $[0,\t_{m,k})$, where $\t_{m,k}:=\inf\{t>0:  X^{m,k}_t<0\}$.
For notational simplicity, we again fix $m$ and $k$ and denote $b=b^{m,k}$ and $\si=\si^{m,k}$, so that (\ref{SDE4}) now reads:
\bea
\label{SDE5}
\left\{\ba{lll}
\dis X_t=x+\int_0^tb(s, X_s, W_s)ds+\int_0^t\si(s, X_s,W_s)dB_s-Q_t,\ms\\
W_t= t-\si_{N_t},
\ea\right.
\eea
Recalling from (\ref{GaXi}) and (\ref{approxstratb}) that the function $b=b^{m,k}$ is discontinuous, but has a linear growth:
\bea
\label{lineargr}
|b(t,x,w)|\leq C(1+|x|),\qq (t,x,w)\in [0,T]\times\hR\times[0,T].
\eea
for some constant $C>0$ depending only on the coefficients but independent of $m,k$. In what follows we shall allow such generic  constant to vary from line to line.

The scheme for constructing the strong solution for (\ref{SDE5}) goes as follows (see, e.g.,  \cite{GP, NO} or \cite{BM}). For any $N>0$ define $b_{N}(t,x,w)=b(t,x\wedge N\vee(-N),w)$.  Then (\ref{lineargr}) implies that $b_{N}$ is a bounded measurable function.
Let $\rho$ be a smooth mollifier with compact support in $\hR$ such that $\int_{\hR}\rho(z)dz=1$. For $n=1,2....$, define
$$b_{N,j}(t,x,w)=j\int b_{N}(t,z,w)\rho (j(x-z))dz,$$
then  $b_{N,j}$'s  are smooth functions, having the same bound $N$, and satisfying  the linear growth condition (\ref{lineargr}) with the same constant $C>0$, and $b_{N, j}\to b_N$ almost everywhere on $[0,T]\times\hR\times[0,T]$, as $j\to\infty$.

Next, for any $K\in\hN$ and $j\leq K$ we define
$\tilde{b}_{N,j,K}\dfnn \bigwedge_{k=j}^{K}b_{N,j}$ and $\tilde{b}_{N,j}\dfnn\bigwedge_{k=j}^{\infty}b_{N,j}$, where $a\wedge b=\min\{a,b\}$. Then clearly, each $\tilde{b}_{N,j,K}$ is continuous, and uniformly Lipschitz in $x$, uniformly in $(t,w)$. Furthernore, for almost all $x$, for any $(t,w)$, it holds that
\bea
\label{bconv}
\tilde{b}_{N,j,K}\downarrow \tilde{b}_{N,j},\  {\rm as } \  K\rightarrow \infty \q \mbox{\rm and}\q
\tilde{b}_{N,j}\uparrow b_{N},\ {\rm as}\  j\rightarrow \infty.
\eea
 Now let us fix $N$, $j$, and $K$, and consider the SDE:
\bea
\label{SDEn}
\left\{\ba{lll}
\dis dY_{t}=\tilde{b}_{N,j,K}(t,Y_{t},W_{t})dt+\si(t,Y_{t},W_{t})dB_{t}, \q Y_0=x; \\
W_t =t-\si_{N_t}, \qq\qq t\ge0.
\ea\right.
\eea
Clearly, (\ref{SDEn}) has a unique strong solution, denote it by $\tilde{Y}^{N,j,K}$. By the standard comparison theorem, we see
that $\{\tilde{Y}^{N,j,K}\}$ is decreasing with $K$, and thus we can define $\tilde{Y}^{N,j}_t\dfnn \lim_{K\to
\infty} \tilde Y^{N,j,K}_t$, $t\in[0,T]$, $\hP$-a.s. Since $\tilde{b}^{N,j}$'s and $\si$ are bounded, one can easily check
that $\tilde{Y}^{N, j}_t<\infty$, $\hP$-a.s.
We shall argue that the limiting process $\tilde{Y}^{N,j}$ solves the SDE:
\bea
\label{SDE6}
\left\{\ba{lll}
\dis dY_{t}=\tilde{b}_{N,j}(t,Y_{t},W_{t})dt+\si(t,Y_{t},W_t)dB_{t}, \q Y_0=x;\\
W_t =t-\si_{N_t},
\ea\right.\qq t\ge0.
\eea
To see this, we first prove the following crucial lemma.
\begin{lem}
\label{Krylov}
Suppose that Assumptions \ref{assump0} and \ref{assump2} are in force. Assume also that
$\{\hat b_{K}\}_{n=1}^{\infty}$ are measurable functions defined on $[0,T]\times\hR\times[0,T]$,  bounded uniformly in $K$, and
there exists a measurable function $\hat b$ such that
$$\lim_{K\to \infty}\hat b_{K}(s,x,w)=\hat b(s,x,w)\ \ \ {\rm for} \ \ a.e. ~(s,x,w)\in [0,T]\times \hR \times [0,T]. $$
Suppose that for each $K$,  $(\hat Y^{K}, W)$ is a strong solution of (\ref{SDEn})  with drift being replaced by $\hat b_{K}$, and
that there exists $\hat Y$ such that for
every $t\in[0,T]$, $\lim_{K\to \infty}\hat Y_{t}^{K}=\hat Y_{t}$, $\hP$-a.s. Then, it holds that
\bea
\label{bconv}
\lim_{K\to \infty}\hE\Big[\int_{0}^{T}|\hat{b}_{K}(t,\hat Y_{t}^{K},W_{t})-\hat{b}(t,\hat Y_{t},W_{t})|ds\Big]=0.
\eea
  \end{lem}

{\it Proof.} Since the proof follows the idea of \cite{NO} or \cite{BM} closely, we only give a sketch for completeness.

 First note that since  each $(\hat Y^K, W)$ is a strong solution to (\ref{SDEn}), we can apply Lemma \ref{krylov1}
and obtain the estimate (\ref{krylovI}) to each $(\hat Y^K, W)$. Note that $\hat b_K$'s are bounded, uniformly in $K$,  we see that
constant $G$ is independent of $K$. Assuming now that the function $g$  in (\ref{krylovI}) is bounded and continuous, a
simple application of Bounded Convergence Theorem then shows that the estimate (\ref{krylovI}) holds for $(\hat Y, W)$. A further
Monotone Class argument then shows that 
the estimate (\ref{krylovI}) hold for $(\hat Y, W)$ actually for any bounded and measurable function $g$.

To prove (\ref{bconv}) we first write
 \bea
\label{In}
\hE\Big[\int_{0}^{T}|\hat b_{K}(t,\hat Y_{t}^{K},W_{t})-\hat b(t,\hat Y_{t},W_{t})|ds\Big]\leq I^K_{1}+I^K_2,
\eea
where
\beaa
\label{In12}
I^K_{1}:= \sup_{k}\hE\Big[\int_{0}^{T}|\hat b_{k}(t,\hat Y_{t}^{K},W_{t})-\hat b_{k}(t,\hat Y_{t},W_{t})|dt\Big], ~
I^K_{2}:= \hE\Big[\int_{0}^{T}|\hat b_{K}(t,\hat Y_{t},W_{t})-\hat b(t,\hat Y_{t},W_{t})|dt\Big].
\eeaa
Let $\kappa:\hR\rightarrow [0,1]$ be a smooth truncation function satisfying
$\kappa(z)=0$ for $|z|\geq 1$ and $\kappa(0)=1$. Then by Bounded Convergence Theorem one has
\bea
\label{kappato0}
\lim_{R\to\infty}\hE\Big[\int_0^T(1-\kappa(\hat Y_{t}/R))dt\Big]=0.
\eea
Now we fix $R>0$, and denote $\D \hat b_K:=\hat b_K-\hat b$. Since both $\hat b_{K}$ and $\hat b$ are bounded and continuous,
we apply Lemma \ref{krylov1} with $g=\D\hat b_K$ to $(\hat Y, W)$ to get
 \bea
 \label{I2n}
I^K_2&=& \hE\Big[\int_{0}^{T} \kappa(\hat Y_{t}/R)|\D\hat b_{K}(t,\hat Y_{t},W_{t})|dt\Big]+ \hE\Big[\int_{0}^{T} (1-\kappa(\hat Y_{t}/R))|\D\hat b_{K}(t,\hat Y_{t},W_{t})|dt\Big]\nonumber\\
&\leq & G\Big(\int_{0}^{T}\int_{-R}^{R}\int_{0}^{t}|\D\hat b_K(t,x,w)|^{2}dwdxdt\Big)^{1/2}+2C\hE\int_{0}^{T}(1-\kappa(X_{t}/R))dt.
\eea
Since $\lim_{K\to\infty}\D\hat b_K= 0$, first letting $K\to \infty$ and then letting $R\rightarrow \infty$ we get $\lim_{K\to\infty}I^K_{2}=0$.

To show $I^K_1 \to 0$, we first note that by (\ref{kappato0}), for any $\e>0$,
there exists $R_0$ such that
\bea
\label{Kappalee}
\hE\Big[\int_{0}^{T}|1-\kappa(\hat Y_{t}/R_0)|dt\Big]<\varepsilon.
\eea
Since $\hat b_{K}\to \hat b$ a.e., as $K\to\infty$, and all $\{\hat b_k\}$'s and $\hat b$ are bounded, it is clear that
$\hat b_{K}\to \hat b$ in $L^2_{T,R_0}:= L^{2}([0,T]\times [-R_0,R_0]\times [0,T])$, hence $\{\hat b_{K}, \hat b\}_{K\ge1}$ is a compact set in  $L^2_{T,R_0}$. Thus, for any $\e>0$, we can find finitely many bounded smooth functions $\{H_{j}\}_{j=1}^L$, such that
for each $k$, there is a $H_{i_k}$ so that
\bea
\label{bkHi}
\Big(\int_{0}^{T}\int_{-R_0}^{R_0}\int_{0}^{t}|\hat b_k(t,x,w)-H_{i_k}(t,x,w)|^{2}dwdxdt\Big)^{1/2}<\e.
\eea
Now, we write
\beaa
I^K_1=\sup_{k}\hE\int_{0}^{T}|\hat b_{k}(t,\hat Y_{t}^K,W_{t})-\hat b_{k}(t,\hat Y_{t},W_{t})|dt\leq \sup_{k}I_{1}(K,k)+I_{2}(K)+\sup_{k}I_{3}(k),
\eeaa
where
\beaa
\left\{\ba{lll}
\dis I_{1}(K,k)=\hE\Big[\int_{0}^{T}|\hat b_{k}(t,\hat Y^K_{t},W_{t})-H_{i_k}(t,\hat Y_{t}^K,W_{t})|dt\Big];\\
\dis I_{2}(K)=\sum_{j=1}^{L}\hE\Big[\int_{0}^{T}|H_{j}(t,\hat Y^K_{t},W_{t})-H_{j}(t,\hat Y_{t},W_{t})|dt\Big];\\
\dis I_{3}(k)=\hE\Big[\int_{0}^{T}|\hat b_{k}(t,\hat Y_{t},W_{t})-H_{i_k}(t,\hat Y_{t},W_{t})|dt\Big].
\ea\right.
\eeaa

By Bounded Convergence Theorem, we have $\lim_{K\rightarrow\infty}I_{2}(K)=0$. Next, similar to (\ref{I2n}) we apply the estimate (\ref{krylovI}) with $\b\gamma=2$, along with (\ref{Kappalee}) and (\ref{bkHi}), to get, for each $k$,
\beaa
I_{1}(K,k)
&\leq & G\Big( \int_{0}^{T}\int_{-R_0}^{R_0}\int_{0}^{t}|[\hat b_{k}-H_{i_k}](t,x,w)|^{2}dwdxdt|\Big)^{1/2}+C_{1} \hE\Big[\int_{0}^{T}(1-\kappa (X_{t}^{n}/R_0))dt\Big]\\
&\le& G\e+C_1\e,
\eeaa
where $G$ is defined by (\ref{G}) with $\b\g=2$, and $C_{1}$ is a constant depending on $C$ and
$\max_{1\leq i \leq L}\|H_{i}\|_{\infty}$. Consequently, we have
$\lim_{K\rightarrow \infty}\sup_{k}I_{1}(K,k)\leq (G+C_{1})\e$.
Similarly, we can show that $\sup_{k}I_{3}(k)\leq (G+C_{1})\e$. Letting $\e\to 0$ we obtain
 $\lim_{K\rightarrow \infty}I^K_1=0$. The proof is now complete.
 \qed

Let us fix $N, j$ and denote $\hat b_K=\tilde b_{N,j,K}$, $\hat Y^K=\tilde Y^{N, j, K}$, $K\in \hN$, and $\hat b=\tilde b_{N,j}$,
$\hat Y=\tilde Y^{N,j}$. Then Lemma \ref{Krylov} shows that, possibly along a subsequence and may assume its own, we have
\bea
\label{limb}
\lim_{K\to\infty} \int_0^t \tilde b_{N, j, K}(s, \tilde Y^{N, j, K}_s, W_s)ds = \int_0^t\tilde b_{N,j}(s, \tilde Y^{N, j}_s, W_s)ds,\q t\in[0,T], \q\hP\as.
\eea
 On the other hand,  since $\si$ is bounded and continuous, by bounded convergence theorem it is easy to see that
$\lim_{K\to\infty} \hE\big[\big|\int_0^T[\si(s, \tilde{Y}^{N, j, K}_s, W_s)-\si(s, \tilde{Y}^{N, j}_s, W_s)dB_s\big|^2\big]=0,
$
hence, modulo a subsequence we have
\bea
\label{limsi}
\lim_{K\to\infty}\int_0^t\si(s, \tilde{Y}^{N, j, K}_s, W_s)dB_s= \int_0^t\si(s, \tilde{Y}^{N, j}_s, W_s)dB_s, \q t\in[0,T], \q \hP\as,
\eea
as well. (\ref{limb}) and (\ref{limsi}), together with the facts that $\tilde Y^{N,j,K}$ solves SDE (\ref{SDEn}) and $\tilde Y^{N, j, K}\da \tilde Y^{N,j}$ show that $\tilde Y^{N,j}$ solves the
SDE (\ref{SDE6}).

Next, since $\tilde{Y}^{N,j,K}\leq \tilde{Y}^{N,i,K}$, for $j\leq i \leq K$,  we see that $\tilde{Y}^{N,j}$
increases as $j$ increases, thus $\tilde{Y}^{N,j}_t\ua Y^N_t$, $t\in[0,T]$, $\hP$-almost surely, where $Y^{N}$ is some process
with $Y^N_t<\infty$, $t\in[0,T]$, $\hP\as$ The same argument as before, using Lemma \ref{Krylov} with $\hat b_j=b_{N,j}$,
$\hat b=b_N$, and $\hat Y^j=Y^{N,j}$, we can show that $Y^{N}$ solves the SDE:
\bea
\left\{\ba{lll}
\dis dY_{t}={b}_{N}(t,Y_{t},W_{t})dt+\si(t,Y_{t},W_{t})dB_{t}, \q Y_0=x;\\
W_t =t-\si_{N_t},
\ea\right. \qq t\in[0,T].
\eea
Moreover, we can show, as \cite{BM}, that $Y^N$ is pathwisely unique. Let us now define $\tau_{N}=\inf\{t: |Y_{t}^{N}|\geq N\}
\wedge T$. Then on the interval $[0,\tau_{N}]$, $b_N(t, Y^N_t, W)=b(t, Y^N_t,W)$, thus $Y^{N}$ is a unique strong solution to the SDE
\bea
\label{SDE7}
\left\{\ba{lll}
\dis dY_{t}= b(t,Y_{t},W_{t})dt+\si(t,Y_{t},W_{t})dB_{t}, \q Y_0=x; \\
W_t =t-\si_{N_t},
\ea\right. \qq t\in [0,\t_N].
\eea
Now observe that if $N_1>N_2$, we have $\t_{N_1}\ge \t_{N_2}$. Thus by uniqueness we have $Y^{N_2}_t=Y^{N_1}_t$ on the
interval $[0,\t_{N_2}]$. We can now define a process $Y$ such that $Y_t=Y^N_t$, $t\in[0,\t_N]$.  Then $Y$ is well-defined on
the interval $[0,\t)$, where  $\tau=\lim_{N\ua \infty}\tau_{N}$. Since $b$ is of linear growth, and $\si$ is bounded, it is not hard
to show that $\hE[\sup_{t\in[0,T]}|Y^N_t|^2]<\infty$, which implies that $\hP\{|Y_t|<\infty, t\in[0,\t)\}=1$, and hence $\t=T$, $\hP$-a.s.
In other words, $Y$ is a unique strong solution to (\ref{SDE7}) on $[0,T]$.

We  can now prove the main result of this section.
 \begin{thm}
Assume that Assumptions \ref{assump0} and \ref{assump2} are in force. Then, for each $k>0$, the closed loop system (\ref{closeloop}) possesses a unique strong solution $(X^k,W)$ on the random interval $[0,\t_k)$, where $\tau_k=\inf\{t>0: X^k< 0\}\wedge T$.
\end{thm}

{\it Proof.} We begin by recalling the SDE (\ref{SDE4}). Without loss of generality
we consider only the case $s=0$, that is, we write  SDE (\ref{SDE4})  as
\bea
\label{SDE8}
\left\{\ba{lll}
\dis dX_{t}= b^{m,k}(t,X_{t},W_{t})dt+\si^{m,k}(t,X_{t},W_{t})dB_{t}-dQ_{t}, \q X_0=x;\ms\\
W_t =t-\si_{N_t}
\ea\right. \q t\in[0,T].
\eea

We shall follow  the same argument as that in Proposition \ref{Xmk} to construct the strong solution on the canonical space
$(\O^1,\cF^1,\hP^0;\hF^1)$ defined by (\ref{WP}). For any $\o=(\o^1,\o^2)\in\O$, we write the coordinate processes as $B_t(\o)\dfnn \o^1(t)$, $ Q_t(\o)\dfnn \o^2(t)$, $(t,\o)\times [0,T]\times\O$. Assuming that  the process $Q_t(\o)=\o^2(t)$ jumps  at $0<\sigma_{1}(\o^2)< \cdots<\sigma_{N_T(\o^2)}(\o^2)<T$, where $N_t(\o^2)$ denotes the number of jumps of $Q$ up to time $t$, we  define  $W_t(\o)=t-\si_{N_t(\o^2)}(\o^2)$, $t\ge 0$.

Now for $\hP^Q$-a.s. $\o^2\in \O^2$ we define $\tilde b^{m,k,\o^2}$ and $\tilde \si^{m,k,\o^2}$ by (\ref{bOx}), respectively,
and consider  the SDE on the space $(\O^1,\cF^1,\hP^0;\hF^1)$:
\bea
\label{SDE-1}
\tilde dX_{t}= b^{\o^2, m,k}(t,\tilde X_t)ds+\tilde \si^{\o^2, m,k}(t,\tilde X_{t})dB_t, \q X_0=x;\qq   t\in[0,T],
 \eea
Clearly, this equation is the same as (\ref{SDE7}), and we have shown that it has a unique strong solution on $(\O^1,\cF^1,\hP^0;\hF^1)$, denote it by
$\tilde X^{m,k,\o^2}_t:=\tilde X^{m,k}_t(\cd, \o^2)$ on $(\O^1,\cF^1,\hP^0;\hF^1)$, for $\hP^Q$-a.s. $\o^2\in\O^2$. We then
 define $X^{m,k}:=\tilde X^{m,k}-Q$, and $W_t=t-\si_{N_t}$, then $(X^{m,k}, W)$ is the unique strong solution to (\ref{SDE8}).

To complete the proof, let us define, for fixed $k$,  $\t_{m,k}:=\inf\{t>0, X^{m,k}_t\notin [\frac1m, m]\}\wedge T$.
Again, observe that $b^{m,k}(t, X^{m,k}_t, W)=b^k(t, X^{m,k}_t,W)$ and $\si^{m,k}(t, X^{m,k}_t, W)=\si^k(t, X^{m,k}_t,W)$.
Thus $(X^{m,k},W)$ is the unique strong solution of (\ref{closeloop}) on $[0, \t_m]$. Furthermore, note that if $m_1>m_2$, then
 $\t_{m_1,k}\ge \t_{m_2,k}$. Thus by uniqueness we have $X^{m_2,k}_t=X^{m_1,k}_t$ on the
interval $[0,\t_{m_2}]$. Thus the process $X^k$ defined by $X^k_t=X^{m,k}_t$, $t\in[0,\t_{m,k}]$ is well-defined, and with the linear growth
of $b^k$ and $\si^k$, we see that $\hE[\sup_{t\in[0,T]}|X^{m,k}_t|^2]<\infty$.  We can then conclude that $X^k$ is the unique strong
solution of SDE (\ref{closeloop}) on the interval $[0,\t_k)$, where  $\tau_k=\lim_{m\ua \infty}\tau_{m,k}=\inf\{t>0: X^k<0\}\wedge T$.
\qed

\section{Verification of the $\e$-Optimality}

Having proved the well-posedness of the closeloop system (\ref{closeloop}), we shall now verify that the strategy defined by (\ref{gkak}) is indeed $\e$-optimal. That is, we need to verify that it does produce the cost functional $V^{n_k, \d_k}$ as desired. We should note
that the auxiliary PDE (\ref{HJB2}) actually does not corresponding to any variation of the original control problem (\ref{Xsxw})--(\ref{V1}), the verification is not automatic.

Recall that our $\e$-optimal strategy is based on the
approximating solution $V^{n, \d}$, guaranteed by Theorem \ref{convg}. More precisely,
let $V^k:=V^{n_{k},\d_{k}}\in \hC^{2}_{loc}([0,T]\times\hR)$ be the solutions of (\ref{HJB2}) as those in Thoerem \ref{convg},
such that
\bea
\label{Vk}
\|{V}^{n_k,\d_k}-V\|_{L^{\infty}(D)}<\e_k \searrow 0, \qq \mbox{as $k\to\infty$.}
\eea

Now let us define $\hat V^k(s,x,w)=V^k(s,x,w)\1_{D}(s,x,w)$. Then $\hat V^k\in \hC^{1, 2,1}(D)$, and
it follows from (\ref{Vk}) that $\|\hat{V}^{k}-V\|_{L^{\infty}(D)}\rightarrow 0$, as $k\to\infty$. Furthermore, by the construction of ${V}^{k}$, we see
that ${V}^{k}_{x+}(s,-\delta,w)> 1$, and hence $\hat{V}^{k}_{x+}(s,0,w)={V}^{k}_{x+}(s,0,w)> 1$ for $k$ large enough.
We should note that $\hat V^k_x(s,x,w)=V_{x}^{k}(s,x,w)>0$ for $(s,x,w)\in D$ always holds.

We now recall the strategy $\pi^k=(\g^k,a^k)$ defined by (\ref{gkak}) and denote $X^k$ be the corresponding strong solution to (\ref{Xsxw}), which exists on $[0, \t^k)$, where $\t^k:=\inf\{t>0: X^k_t\notin [0,\infty)\}$. It is useful to remember that
$\pi^k$ is actually the maximizer of the Hamiltonian (\ref{H0}), namely, it holds that
\bea
\label{gkmax}
%
\g^{k}_{t}=\argmax_{\gamma\in [0,1]}\Big[\frac12\si^{2}\g^{2} (X^{k}_{t})^{2}\hat V_{xx}^{k}(t,X^{k}_{t},W_{t})+(\mu-r) \gamma X^{k}_{t}
\hat V_{x}^{k}(t,X^{k}_{t},W_{t})\Big].
\eea
In the rest of the section we shall consider, for $s\in[0,T]$,  the closeloop system (\ref{closeloop}) on the interval $[s, T]$, and write it as:
\bea
\label{Xsxw1}
\left\{\ba{lll}
dX_t=b^k(t, X_t)dt+\si^k(t, X_t)dB_t-dQ^{s,w}_t; \qq X_s=x;\\
W_t=w+(t-s)-(\si_{N_t}-\si_{N_s}),
\ea\right. \q t\in[s, T].
\eea
where
$b^k(t,x)=(p-a^k_t)+[r+(\m-r)\g^k_t]x$; $\si^k(t,x)=\g^k_tx$, and $\pi^k=(\g^k, a^k)$ is the aforementioned approximating strategy.
We denote the solution  by $X^{k}=X^{k, s,x}$ and $W=W^{s,w}$ when the context is clear. For given $(s,x,w)\in D$ we define
 $\t^k_s:=\inf\{t>s: X^k_t\notin [0,\infty)\}$,
 and denote $\hE_{sxw}[\,\cd\,]:=\hE[\,\cd\,|X^k_s=x, W_s=w]$.

To show that the strategy $\pi^k=(\g^k, a^k)$ does satisfy the $\e$-optimality we shall argue that $J(s,x,w; \pi^k)$ satisfies, for all $(s,x,w)\in D$, that
\bea
\label{epop}
J(s,x,w;\pi^k)\to V(s,x,w), \qq \mbox{\rm as $k\to \infty$.}
\eea
But note that $J(s,x,w;\pi^k)=\hE_{sxw}\big[\int_0^{\t^k_s\wedge T} e^{-c(t-s)}a^k_tdt\big]$, and $\lim_{k\to\infty}\|V^k-V\|_{L^\infty(D)}= 0$,
the following theorem would be sufficient.
\begin{thm}
Assume that Assumptions \ref{assump0}
and \ref{assump2} are in force. Then it holds that
 \bea
 \label{92901}
\lim_{k \rightarrow \infty }\Big|\hE_{sxw}\Big[\int_{s}^{{{\tau}}^{k}_s\wedge T}e^{-c(t-s)}a^{k}_{t}dt-\hat{V}^{k}(s,x,w)\Big]\Big|=0,  \qq \mbox{\rm uniformly in $(s,x,w)\in  D$.}
\eea
\end{thm}
{\it Proof.} The proof is straightforward. We first apply It\^o's formula from $s$ to $\t^k_s\wedge T$ to the process $e^{-c(t-s)}\hat{V}^{k}(t,X^{k}_t,W_t)$ to get
\beaa
&&e^{-c(\tau^{k}_s\wedge T-s)}\hat{V}^{k}(\tau^{k}_s\wedge T,X^{\k}_{\tau^{k}_s\wedge T},W_{\tau^{k}_s\wedge T} )\nonumber \\
&=&\hat{V}^{k}(s,x,w)+ \int_{s}^{\tau^{k}_s\wedge T}e^{-c(t-s)}\Big[-c\hat{V}^{k}+\hat{V}^{k}_{t}+\hat{V}^{k}_{w}+[(p-a^{k}_{t})+(r+(\mu-r)\g^k_{t}) X^{k}_{t}] \hat{V}^{k}_{x} \nonumber\\
&&+ \frac12\si^{2}(\g^{k}_{t})^{2} (X^{k}_{t})^{2}\hat{V}^{k}_{xx}\Big](t,X^{k}_{t},W_{t})dt+\si\int_{s}^{\tau^{k}_s\wedge T}e^{-c(t-s)} \gamma^{k}_{t}X_{t}^{k}dB_{t}\nonumber\\
&& +\sum_{s\leq t\leq \tau^{k}_s\wedge T}e^{-c(t-s)}(\hat{V}^{k}(t,X^{k}_{t},W_{t})-\hat{V}^{k}(t,X^{k}_{t-},W_{t-})).
\eeaa
Taking the expectation on both sides above yields
\beaa
&&\hE\Big[e^{-c(\tau^{k}_s\wedge T-s)}\hat{V}^{k}(\tau^{k}_s\wedge T,X^{\k}_{\tau^{k}_s\wedge T},W_{\tau^{k}_s\wedge T} )\Big]\nonumber \\
&=&\hat{V}^{k}(s,x,w)+ \hE\Big[\int_{s}^{\tau^{k}_s\wedge T}e^{-c(t-s)}\Big[-c\hat V^k+\hat{V}^{k}_t+\hat{V}^{k}_{w}+[(p-a^{k}_{t})+(r+(\mu-r)\g^k_{t}) X^{k}_{t}] \hat{V}^{k}_{x} \nonumber\\
&&+ \frac12\si^{2}(\g^{k}_{t})^{2} (X^{k}_{t})^{2}\hat{V}^{k}_{xx}\Big](t,X^{k}_{t},W_{t})dt\Big]\nonumber\\
&&+\hE\Big\{\int_{s}^{\tau^{k}_s\wedge T}e^{-c(t-s)}\frac{f(W_{t})}{\overline{F}(W_{t})}\Big[\int_{0}^{X^{k}_{t}} \hat{V}^{k}(t,X^{k}_{t}-u,0)g(u)du-
\hat{V}^{k}(t,X^{k}_{t},W_{t})\Big]dt\Big\}.
\eeaa
Since $\hat{V}^{k}(s,x,w)$ satisfies the HJB equation (\ref{HJB3}), and $\pi^k=(\g^k,a^k)$ is the maximizer in terms of $\hat V^k$,
a simple calculation shows that (suppressing variables)
\beaa
&& \neg-\neg c\hat{V}^{k}+\hat{V}^{k}_{t}+\hat{V}^{k}_{w}+[(p-a^{k}_{t})+(r+(\mu-r)\g^k_{t}) X^{k}_{t}] \hat{V}^{k}_{x}+ \frac12\si^{2}(\g^{k}_{t})^{2} (X^{k}_{t})^{2}\hat{V}^{k}_{xx}-\frac{f(W_{t})}{\overline{F}(W_{t})}\hat{V}^{k}\\
&=&-a^{k}_{t}-\frac{f(W_{t})}{\overline{F}(W_{t})}\int_{0}^{X^{k}_{t}+\delta_{k}} V^{k}(t,X^{k}_{t}-u,-\delta_{k})g(u)du-\frac{\e_{k}}{2}\hat{V}^{k}_{xx} -\frac{\e_{k}}{2}\hat{V}^{k}_{ww}.
\eeaa
Then we have
\beaa
&&\hE\Big[e^{-c(\tau^{k}_s\wedge T-s)}\hat{V}^{k}(\tau^{k}_s\wedge T,X^{\k}_{\tau^{k}_s\wedge T},W_{\tau^{k}_s\wedge T} )\Big]-\hat{V}^{k}(s,x,w)+\hE\Big[\int_{s}^{\tau^{k}_s\wedge T}e^{-c(t-s)}a^{k}_{t}dt\Big]\nonumber \\
&=&
\hE\Big\{\int_{s}^{\tau^{k}_s\wedge T}e^{-c(t-s)}\frac{f(W_{t})}{\overline{F}(W_{t})}\times\nonumber\\
&&\Big[\int_{0}^{X^{k}_{t}} [ \hat{V}^{k}(t,X^{k}_{t}-u,0)-V^{k}(t,X^{k}_{t}-u,-\delta_{k})]g(u)du-\int_{X^{k}_{t}}^{X^{k}_{t}+\delta_{k}}V^{k}(t,X^{k}_{t}-u,-\delta_{k})g(u)du\Big]\Big\}\\
&&- \frac{\e_{k}}{2}\hE\Big[\int_{s}^{\tau^{k}_s\wedge T}e^{-c(t-s)}\hat{V}^{k}_{xx}(t,X^{k}_{t},W_{t}) dt \Big]- \frac{\e_{k}}{2}\hE\Big[\int_{s}^{\tau^{k}_s\wedge T}e^{-c(t-s)}\hat{V}^{k}_{ww}(t,X^{k}_{t},W_{t}) dt\Big]\\
&\leq& C \delta_{k} - \frac{\e_{k}}{2}\hE\Big[\int_{s}^{\tau^{k}_s\wedge T}e^{-c(t-s)}\hat{V}^{k}_{xx}(t,X^{k}_{t},W_{t}) dt\Big] -
\frac{\e_{k}}{2}\hE\Big[\int_{s}^{\tau^{k}_s\wedge T}e^{-c(t-s)}\hat{V}^{k}_{ww}(t,X^{k}_{t},W_{t}) dt\Big].
\eeaa
Now, letting $k\to \infty$ and noting that $\d_k, \e_k\to 0$, we see that (\ref{92901}) follows from the fact that
 \beaa
 \label{92902}
\lim_{k\to\infty }\hE_{sxw}\big[\hat{V}^{k}(\tau^{k}_s\wedge T,X^{k}_{\tau^{k}_s\wedge T},W_{\tau^{k}_s\wedge T} )\big]
&=&\lim_{k\to\infty }\hE_{sxw}\big[\1_{\{\tau^{k}_s\geq T\}}\hat{V}^{k}( T,X^{k}_{T},W_{T} )\big]\\
&=&\lim_{k\to\infty }\hE_{sxw}\big[\1_{\{\tau^{k}_s\geq T\}}V( T,X^{k}_{T},W_{T} )\big]=0.
\eeaa
 This proves the theorem.
\qed

%

\section{Appendix}
\setcounter{equation}{0}


[{\it Proof of Lemma \ref{classL1}.}]
We shall prove only the sub-solution case, the super-solution case is similar to \cite{GMS}.
First recall the distance function
$d(x; D):=\inf_{y\in D} |x-y|$, for $x\in\hR^m$, and $D\subset \hR^m$.

%
%
%

Now for any $\th:=(s,x,w)\in(0,T+\delta)\times \hR^{2}$, we define $d_{\sD_{\delta}}(\th):=d(\th; \sD_{\delta}^{c})$. Define the function
\bea
\label{psi}
\psi(\th):=-kd_{\sD_{\delta}}(\th), \qq \th\in (0,T+\d)\times \hR^2,
\eea
where, recalling the constant $b$ defined by (\ref{b}), and for $s\in[0,T]$ and $\th=(s,x,w)$, denoting $\G^\d_s:= [0,T]\times [-\delta, T+3\delta]\times [0,s]$, 
\bea
\label{k}
0<k\leq \min\Big\{b-1,~M-K_{2}, ~\inf_{\th\in\G^\d_s}\Psi_{x}(\th), ~~\frac{M-K_{2}}{ \sup_{w\in[0,T+1]}\big|c+\frac{f(w)}{\bar{F}(w)}-r\big|(T+4\delta)}\Big\}.
\eea

 We shall argue that $\psi+\Psi$ will be a viscosity subsolution
of class ($\Psi$) in the sense of Definition \ref{classP}.
To see this, let us first observe that by definition of $\sD^c_{\d}$, one can easily check that
\bea
\label{disD}
 d_{\sD_\d}(s,x,w)=(x+\d)\wedge (w+\delta) \wedge \frac{\sqrt{2}}{2}(s+\delta-w)\wedge (T+\delta-s)\wedge s, \q
 \mbox{\rm for  $(s,x,w)\in \sD_{\delta}$.}
\eea
So for any $\th:=(s,x,w)\in \sD_\d$, we shall consider the following cases:

\ms

{\bf Case 1.}  $d_{\sD_{\delta}}(\th)=x+\delta<(w+\delta) \wedge \frac{\sqrt{2}}{2}(s+\delta-w)\wedge (T+\delta-s)\wedge s$. In this case,
$x<T+3\delta$.
Then by definition
(\ref{psi}) and the constraint (\ref{k}) we have
\bea
\label{case1}
&& (\psi+\Psi)_{t}(\th)+\sup_{\g\in [0,1]\atop a\in [0,M]}H^{n}(\th,\psi+\Psi,\nabla (\psi+\Psi)_x,(\psi+\Psi)_{xx}, (\psi+\Psi)_{ww},I^{\delta}[\psi+\Psi],\g, a)\nonumber\\
&\geq& -K_{2}-k(r(x+\delta)+p)+\Big(c+\frac{f(w)}{\bar{F}(w)}\Big)k(x+\delta)+M(1+k)\\
&\geq& k\Big[\Big(c+\frac{f(w)}{\bar{F}(w)}-r\Big)(x+\delta)+M-p\Big]+M-K_{2} \nonumber\\
&\geq& k\Big[\Big(c+\frac{f(w)}{\bar{F}(w)}-r\Big)(x+\delta)\Big]+M-K_{2}\geq 0. \nonumber
\eea
where the last inequality follows from (\ref{k}).

\ms
{\bf Case 2.} $d_{\sD_{\delta}}(\th)=w+\delta<(x+\delta)\wedge \frac{\sqrt{2}}{2}(s+\delta-w)\wedge (T+\delta-s)\wedge s$. In this case we have
\bea
&&(\psi+\Psi)_{t}(\th)+\sup_{\g\in [0,1]\atop a\in [0,M]}\sH^{n}(\th,\psi+\Psi,\nabla (\psi+\Psi)_x,(\psi+\Psi)_{xx}, (\psi+\Psi)_{ww},I^{\delta}[\psi+\Psi],\g, a)\nonumber\\
&& \qq \geq -K_{2} -k+M+k\Big(c+\frac{f(w)}{\bar{F}(w)}\Big)(w+\delta)\geq -K_{2} -k+M\geq 0,
\eea
again, thanks to (\ref{k}).

\ms
{\bf Case 3.} $d_{\sD_{\delta}}(\th)=\frac{\sqrt{2}}{2}(s+\delta-w)<(x+\delta)\wedge (w+\delta)\wedge (T+\delta-s)\wedge s$.
Similarly, using Assumption \ref{assump1}, we can calculate that
\bea
&& (\psi+\Psi)_{t}(\th)+\sup_{\g\in [0,1]\atop a\in [0,M]}\sH^{n}(\th,\psi+\Psi,\nabla (\psi+\Psi)_x,(\psi+\Psi)_{xx}, (\psi+\Psi)_{ww},I^{\delta}[\psi+\Psi],\g, a)\nonumber\\
&&  \qq \qq \ge-K_{2}+M\geq 0,
\eea

\ms
{\bf Case 4.}  $d_{\sD_{\delta}}(\th)=(T+\delta-s)<(x+\delta)\wedge (w+\delta) \wedge \frac{\sqrt{2}}{2}(s+\delta-w)\wedge s $. Again, we have
\bea
\label{case4}
&&(\psi+\Psi)_{t}(\th)+\sup_{\g\in [0,1]\atop a\in [0,M]}\sH^{n}(\th,\psi+\Psi,\nabla (\psi+\Psi)_x,(\psi+\Psi)_{xx}, (\psi+\Psi)_{ww},I^{\delta}[\psi+\Psi],\g, a)\nonumber\\
&& \qq \qq \ge -K_{2}+M+k\geq 0,
\eea
where the inequality is again due to (\ref{k}).

 \ms
{\bf Case 5.} $d_{\sD_{\delta}}(\th)=s<(x+\delta)\wedge w+\delta\wedge \frac{\sqrt{2}}{2}(s+\delta-w)\wedge (T+\delta-s)$. In this case we have
\bea
\label{case5}
&&(\psi+\Psi)_{t}(\th)+\sup_{\g\in [0,1]\atop a\in [0,M]}\sH^{n}(\th,\psi+\Psi,\nabla (\psi+\Psi)_x,(\psi+\Psi)_{xx}, (\psi+\Psi)_{ww},I^{\delta}[\psi+\Psi],\g, a)\nonumber\\
&&\qq\qq \geq  -K_{2} -k+M\geq 0,
\eea
again, thanks to (\ref{k}).

\ms
Finally, we note that the function $d_{\sD_\d}$ could also  take possible values from the following sets:
\bea
\label{Bs}
\left\{\ba{lllll}
B^1:=\{\th:=(s,x,w)\in \sD_{\delta}: d_{\sD_{\delta}}(\th)=x+\delta=w+\delta\};\nonumber\\
B^2:=\{\th\in \sD_\d: d_{\sD_{\delta}}(\th)=x+\delta=\frac{\sqrt{2}}{2}(s+\delta-w)\}; \nonumber\\
B^3:=\{\th\in \sD_\d: d_{\sD_{\delta}}(\th)=x+\delta=T+\delta-s\};\nonumber\\
B^4:=\{\th\in \sD_\d:  d_{\sD_{\delta}}(\th)=\frac{\sqrt{2}}{2}(s+\delta-w)=T+\delta-s\}; \\
B^5:=\{\th\in \sD_\d:  d_{\sD_{\delta}}(\th)=w+\delta=\frac{\sqrt{2}}{2}(s+\delta-w)\}; \nonumber\\
B^6:=\{\th\in \sD_\d:  d_{\sD_{\delta}}(\th)=w+\delta=T+\delta-s \};\nonumber\\
 B^7:=\{\th\in \sD_\d:  d_{\sD_{\delta}}(\th)=x+\delta=s \};\nonumber\\
B^8:=\{\th\in \sD_\d:  d_{\sD_{\delta}}(\th)=w+\delta=s \};\nonumber\\
B^9:=\{\th\in \sD_\d:  d_{\sD_{\delta}}(\th)=\frac{\sqrt{2}}{2}(s+\delta-w)=s \};\nonumber\\
 B^{10}:=\{\th\in \sD_\d:  d_{\sD_{\delta}}(\th)=T+\delta-s=s \}.\nonumber
\ea\right.
\eea
Setting ${\sB_{\delta}}:=\cup_{i=1}^{10} B^i$, it is easy to verify that if $\th\in {\sB_{\delta}}$, then one of the following
 must hold:
 \beaa
 \psi_{x}(s,x-,w)<\psi_{x}(s,x+,w), ~ \mbox{\rm  or}~  \psi_{w}(s,x,w-)<\psi_{w}(s,x,w+), ~ \mbox{\rm  or}~   \psi_{s}(s-,x,w)<\psi_{s}(s+,x,w).
 \eeaa
That is, if $\th\in {\sB_{\delta}}$, then the $\nabla \psi$ will have a positive jump at $\th$ in one of the directions, or it is ``convex" at $\th$ in that direction. Therefore one
cannot find any smooth test function $g$ that is above $\psi$ so that
$0=(\psi-g)(\th)$ is a strict maximum
over $\sD_{\delta}$. This, together with (\ref{case1})--(\ref{case5}), shows that $\underline{\psi}:=\psi+\Psi$ is a viscosity subsolution  to
 (\ref{HJB2}). Furthermore, by  definition of $\psi$ and $\Psi$, it is readily seen that $\underline{\psi}$ is  of class ($\Psi$). This
 proves the lemma.
\qed

\end{document}